\documentclass[12pt]{article}
\usepackage[final]{epsfig}
\usepackage{graphics}
\usepackage{amsmath,amssymb,amsthm,amscd}
\usepackage{latexsym}
\usepackage{epsfig}
\usepackage{amsfonts}
\usepackage{latexsym}
\usepackage{graphicx}
\usepackage{epstopdf}
\usepackage{multicol,multirow}
\makeatletter
\newtheorem*{rep@theorem}{\rep@title}
\newcommand{\newreptheorem}[2]{%
\newenvironment{rep#1}[1]{%
 \def\rep@title{#2 \ref{##1}}%
 \begin{rep@theorem}}%
 {\end{rep@theorem}}}
\makeatother

\newtheorem{lemma}{Lemma}[section]
\newtheorem{proposition}[lemma]{Proposition}
\newtheorem{remark}[lemma]{Remark}

\newtheorem{theorem}[lemma]{Theorem}
\newtheorem{definition}[lemma]{Definition}
\newtheorem{corollary}[lemma]{Corollary}

\newtheorem{proposition-conjecture}[lemma]{Proposition-conjecture}
\newtheorem{problem}[lemma]{Problem}
\newtheorem{thma}{Theorem}

\newreptheorem{thma}{Theorem}

\begin{document}
\newcommand{\eps}{{\varepsilon}}
\newcommand{\proofend}{\hfill$\Box$\bigskip}
\newcommand{\C}{{\mathbb C}}
\newcommand{\Q}{{\mathbb Q}}
\newcommand{\R}{{\mathbb R}}
\newcommand{\Z}{{\mathbb Z}}
\newcommand{\RP}{{\mathbb {RP}}}
\newcommand{\CP}{{\mathbb {CP}}}
\newcommand{\PP}{{\mathbb {P}}}
\newcommand{\ep}{\epsilon}
\newcommand{\G}{{\Gamma}}

\def\bull{\\$\bullet$ \thinspace}

\def\proof{\paragraph{Proof.}}
\def\prooflem{\paragraph{Proof of lemma.}}
\def\proofthm{\paragraph{Proof of theorem.}}
\def\proofprop{\paragraph{Proof of proposition.}}
\def\proofsketch{\paragraph{Proof sketch.}}


\newcommand{\marginnote}[1]
{
}

\newcounter{bk}
\newcommand{\bk}[1]
{\stepcounter{bk}$^{\bf BK\thebk}$%
\footnotetext{\hspace{-3.7mm}$^{\blacksquare\!\blacksquare}$
{\bf BK\thebk:~}#1}}

\newcounter{fs}
\newcommand{\fs}[1]
{\stepcounter{fs}$^{\bf FS\thefs}$%
\footnotetext{\hspace{-3.7mm}$^{\blacksquare\!\blacksquare}$
{\bf FS\thefs:~}#1}}


\title {Integrability of higher pentagram maps}

\author{Boris Khesin\thanks{
School of Mathematics, Institute for Advanced Study, Princeton, NJ 08540, USA and Department of Mathematics,
University of Toronto, Toronto, ON M5S 2E4, Canada;
e-mail: \tt{khesin@math.toronto.edu}
}
\,  and Fedor Soloviev\thanks{
Department of Mathematics,
University of Toronto, Toronto, ON M5S 2E4, Canada;
e-mail: \tt{soloviev@math.toronto.edu}
}
\\
}

\date{}

\maketitle

\begin{abstract}
We define higher pentagram maps on polygons in $\PP^d$ for any dimension $d$, which extend R.~Schwartz's definition of the 2D pentagram map.
We prove their integrability by presenting Lax representations
with a spectral parameter for scale invariant maps.
The corresponding continuous limit of the pentagram map in dimension $d$ is shown to be the
$(2,d+1)$-equation of the KdV hierarchy, generalizing the Boussinesq equation in 2D.
We also study in detail the 3D case, where we prove integrability for both closed and twisted polygons and describe the  spectral curve, first integrals, the corresponding tori and the motion along them,
as well as an invariant symplectic structure.
\end{abstract}

\tableofcontents

\section{Introduction} \label{intro}

The pentagram map was defined by R.~Schwartz in \cite{Schwartz} on plane convex 
polygons considered
modulo projective equivalence. Figure 1 explains the definition: for a polygon $P$ the image under the pentagram map is a new polygon $T(P)$ spanned by the ``shortest" diagonals of $P$. Iterations of this map on classes of projectively equivalent polygons manifest quasiperiodic behaviour, which indicates hidden integrability \cite{Sch08}.

\begin{figure}[hbtp]
\centering
\includegraphics[width=1.8in]{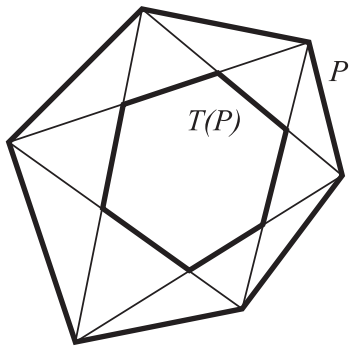}
\caption{\small The image $T(P)$ of a hexagon $P$ under the 2D pentagram map.}
\end{figure}

The integrability was proved in \cite{OST99} for the pentagram map on a larger class of the so called
twisted polygons in 2D, which are piecewise linear curves with a fixed monodromy relating their ends.
Closed polygons correspond to the monodromy given by the identity transformation.
It turned out that there is an invariant Poisson structure
for the pentagram map and it has sufficiently many invariant quantities. Moreover, this map turned out to be related to a variety of mathematical domains, including cluster algebras \cite{FZ, Glick}, frieze patterns,
and integrable systems of mathematical physics: in particular, its continuous limit in 2D is the classical Boussinesq equation \cite{OST99}.
Integrability of the pentagram map for 2D closed polygons was established in \cite{FS11, OST11}, while
a more general framework related to surface networks was presented in \cite{GSTV}.

In this paper we extend the definition of the pentagram map to closed and twisted polygons in spaces of any dimension $d$ and prove its various integrability properties. It is worth mentioning that the problem of finding integrable higher-dimensional generalizations for the pentagram map attracted much attention after the 2D case was treated in \cite{OST99}.\footnote{There seem to be no natural generalization of the pentagram map to polytopes in higher dimension $d\ge 3$. Indeed, the initial polytope should be simple for its diagonal hyperplanes to be well defined. In order to iterate the pentagram map the dual polytope has to be simple as well. Thus iterations could be defined only for  $d$-simplices, which are all projectively equivalent.}
The main difficulty  in higher dimensions is that diagonals of a polygon are generically    skew and do not intersect. One can either confine oneself to special polygons (e.g., corrugated ones, \cite{GSTV}) to retain the intersection property or one has too many possible choices for using hyperplanes as diagonals, where it is difficult to find integrable ones, cf. \cite{Beffa}.

Below, as an analog of the 2D shortest diagonals for a generic polygon in a projective space $\RP^d$
we propose to consider a ``short-diagonal hyperplane" passing through $d$ vertices where every other vertex is taken starting
with a given one. Then a new vertex is constructed as the intersection of $d$ consecutive
diagonal hyperplanes. We repeat this procedure starting with the next vertex of the initial polygon.
The {\it higher} (or {\it $d$-dimensional}) {\it pentagram map} $T$ takes the initial polygon to the one defined by this  set
of new vertices. As before, the obtained polygon is considered modulo projective equivalence
in $\RP^d$.

We also describe  general pentagram maps $T_{p,r}$ in $\RP^d$
enumerated by two integral parameters
$p$ and $r$ by considering $p$-diagonals (i.e., hyperplanes passing through every $p$th vertex of the polygon)
and by taking the intersections of every $r$th hyperplane like that.
There is a curious duality between them: the map $T_{p,r}$ is equal to  $T^{-1}_{r,p}$ modulo a shift in vertex indices.
However, we are mostly interested in the higher pentagram maps, which correspond to $T:=T_{2,1}$ in $\RP^d$.

We start by describing the continuous limit of the higher pentagram map as the evolution in the direction of
the ``envelope" for such a sequence of short-diagonal planes as the number of vertices of the polygon tends to infinity. (More precisely, the envelope here is the curve whose osculating planes are limits of the short-diagonal planes.)

\begin{thma}{\rm ({\bf = Theorems \ref{evolution-anyD}, \ref{thm:kdv}})} \label{thm:cont-lim}
The continuous limit of the higher pentagram map in $\RP^d$ is the $(2,d+1)$-equation in the KdV hierarchy, which is an infinite-dimensional completely integrable system.
\end{thma}

This generalizes the Boussinesq equation as a limit of the pentagram map in $\RP^2$ and
this limit seems to be very robust.
Indeed, the same equation appears for an almost arbitrary choice of diagonal planes.
It also arises when instead of osculating planes one considers other possible definitions of higher pentagram maps (cf. e.g. \cite{Beffa}).

\smallskip

However, the pentagram map in the above definition with short-diagonal hyperplanes  exhibits integrability properties not only in the continuous limit, but as a discrete system as well.
To study them, we define two coordinate systems for twisted polygons  in 3D
(somewhat similar to the ones used in 2D, cf. \cite{OST99}),
and present explicit formulas for the 3D pentagram map using these coordinates (see Theorem \ref{explicit}).

\smallskip

Then we describe the pentagram map as a completely integrable discrete dynamical system by presenting its  Lax form in any dimension and studying
 in detail the 3D case (see Section \ref{sect:alg-geo}). For  algebraic-geometric integrability we complexify the pentagram map.
The corresponding 2D case was  investigated in~\cite{FS11}.

The key ingredient of the algebraic-geometric integrability for a discrete dynamical system is a discrete Lax (or zero curvature) equation with a spectral  parameter, which in our case assumes  the following form:
$$
 L_{i,t+1}(\lambda) = P_{i+1,t}(\lambda) L_{i,t}(\lambda) P_{i,t}^{-1}(\lambda).
$$
 Here the index $t$ represents the discrete time variable, the index $i$ refers to the vertex of an
 $n$-gon, and $\lambda$ is a complex spectral parameter.
 (For the pentagram map in $\CP^d$ the functions $L_{i,t}(\lambda)$ and $P_{i,t}(\lambda)$ are matrix-valued of size $(d+1) \times (d+1)$.)
 The discrete Lax equation arises as a compatibility condition of an over-determined system of equations:
 \begin{equation*}
 \begin{cases}
  L_{i,t}(\lambda) \Psi_{i,t}(\lambda) = \Psi_{i+1,t}(\lambda)\\
  P_{i,t}(\lambda) \Psi_{i,t}(\lambda) = \Psi_{i,t+1}(\lambda),
 \end{cases}
 \end{equation*}
 for an auxiliary function  $\Psi_{i,t}(\lambda)$.

 \begin{remark}
 {\rm
 Recall that for a smooth dynamical system the Lax form is a differential equation of type
$\partial_t L=[P,L]$ on a matrix $L$. Such a form of the equation implies that the evolution of $L$ changes it to a similar matrix, thus preserving its eigenvalues. If the matrix $L$ depends on a parameter, $L=L(\lambda)$, then the corresponding eigenvalues as functions of parameter do not change and in many cases  provide sufficiently many first integrals for complete integrability  of such a system.

Similarly, an analogue of the Lax form for differential operators
of type $\partial_x-L$   is a zero curvature equation $\partial_t
L-\partial_x P=[P,L]\,.$
 This is a compatibility condition which provides the existence of an auxiliary function $\psi=\psi(t,x)$
 satisfying a system  of  equations $\partial_x \psi= L\psi$ and $\partial_t \psi= P\psi\,.$ The above Lax form  and  auxiliary system are discrete versions of the latter.
}
\end{remark}

In our case, the equivalence of formulas for the pentagram map  to the dynamics defined by
the Lax equation implies   complete  algebraic-geometric integrability of the system.
More precisely, the following theorem  summarizes several main results on the 3D pentagram map,
which are obtained by studying its Lax equation.
The dynamics is (generically) defined on the  space ${\mathcal P}_n$ of projectively equivalent twisted $n$-gons in 3D,
which we describe below, and has dimension $3n$, while closed $n$-gons form a submanifold of codimension 15 in it.

Later on we will introduce the notion of spectral data which
consists of a Riemann surface, called a spectral curve, and a
point in the Jacobian (i.e., the complex torus) of this curve, as well as a notion of a spectral map
between the space ${\mathcal P}_n$ and the spectral data.

\begin{thma}{\rm ({\bf = Theorems \ref{spectral-th1}, \ref{time-evol}, \ref{thm:closed}})}
A Zariski open subset of the complexified space ${\mathcal P}_n$ of twisted $n$-gons in 3D is a fibration whose fibres are Zariski open subsets of tori.
These tori are Jacobians of the corresponding spectral curves and are invariant with respect to the space pentagram map.
Their dimension is $3\lfloor n/2\rfloor$ for odd $n$ and $3(n/2)-3$ for even $n$, where
$\lfloor n/2 \rfloor$ is the  integer part of  $n/2$.

The pentagram dynamics on the Jacobians goes along a straight line for odd $n$ and along
a staircase for even $n$ (i.e., the discrete evolution is either a constant shift on a torus,
 or its square is a constant shift).

For closed $n$-gons the tori have dimensions $3\lfloor n/2\rfloor-6$ for odd $n$
and $3(n/2)-9$ for even $n$.
\end{thma}

 \begin{remark}
 {\rm
One also has an explicit description of the the above fibration in terms of coordinates on the space of $n$-gons. We note that the pentagram dynamics understood as a shift on
complex tori does not prevent the corresponding orbits on the
space ${\mathcal P}_n$ from being unbounded. 
The dynamics described above takes place for generic initial
data, i.e., for points on the Jacobians whose orbits do not
intersect certain divisors. Points of generic orbits with irrational shifts can return arbitrarily close to
such divisors. On the other hand, the inverse spectral map is
defined outside of these special divisors and may have poles there.
Therefore the sequences in the space ${\mathcal P}_n$ corresponding to such orbits
may escape to infinity.
}
\end{remark}

\smallskip

It is known that the pentagram map in 2D possesses an invariant Poisson structure \cite{OST99}, which can also be
described by using the Krichever-Phong universal formula \cite{FS11}.
Although we do not present an invariant Poisson structure for the pentagram map in 3D, we
describe its symplectic leaves,  as well as the action-angle coordinates.
More precisely, we present an invariant symplectic structure (i.e., a closed non-degenerate 2-form), and submanifolds
where it is defined (Theorem \ref{thm:rank}).
By analogy with the 2D case, it is natural to suggest that these submanifolds are symplectic leaves of an invariant Poisson structure, and
that the inverse of our symplectic structure coincides with the Poisson structure on the leaves.
An explicit description of this Poisson structure in 3D is still an open problem.

Note that the algebraic-geometric integrability of the pentagram map implies its Arnold--Liouville complete integrability on generic symplectic leaves (in the real case). Namely, the existence of a (pre)symplectic structure coming from the Lax form of the pentagram map (see Section \ref{S:lax-sform}), together with the generic set of first integrals, appearing as coefficients of the corresponding spectral curve, provides sufficiently many integrals in involution. (Note that proving  independence of
first integrals while remaining within the real setting is often more difficult than first proving the algebraic-geometric
integrability,  which  in turn implies their independence in the real case.)

\smallskip

Finally, in Section~\ref{S:higher-lax} we present a Lax form for the pentagram maps in
arbitrary dimension (which implies their complete integrability) assuming their scaling invariance:

\begin{thma}{\rm ({\bf = Theorem \ref{thm:lax_anyD}})}
The scale-invariant pentagram map in $\CP^d$ admits a Lax representation with a spectral parameter.
\end{thma}

The scaling invariance of the pentagram maps is proved for all $d\le 6$, with some numerical evidence for higher values of $d>6$ as well. It would be interesting to  establish it in full generality.
There is a considerable difference between the cases of even and odd dimension $d$, which can be already seen in the analysis of the 2D and 3D cases.

\bigskip

{\bf Acknowledgments}.
We are grateful to M.~Gekhtman and S.~Tabachnikov for useful discussions.
B.K. was partially supported by the Simonyi Fund and an NSERC research grant.


\section{Review of the 2D pentagram map}

In this section we recall the main definitions and results in 2D (see \cite{OST99}), which will be important
for higher-dimensional generalizations below. We formulate the geometric results in the real setting, while the algebraic-geometric ones are presented for the corresponding complexification.

First note that the pentagram map can be extended from closed to twisted polygons.

\begin{definition}
{\rm
Given a projective transformation $M\in PSL(3,\R)$ of the plane $\RP^2$, a {\it twisted $n$-gon} in $\RP^2$  is a map $\phi:\Z\to\RP^2$, such that $\phi(k+n)=M\circ \phi(k)$ for any $k$.
$M$ is called the monodromy of $\phi$. Two twisted $n$-gons are {\it equivalent} if there is a transformation $g\in PSL(3,\R)$ such that $g\circ \phi_1=\phi_2$.
}
\end{definition}

Consider generic $n$-gons, i.e., those that do not have any three consecutive vertices lying on the same line.
Denote by ${\mathcal P}_n$ the space of generic twisted $n$-gons considered up to $PSL(3,\R)$ transformations.
The dimension of ${\mathcal P}_n$ is $2n$. Indeed, a twisted $n$-gon
depends on $2n$ variables representing coordinates of vertices $v_k:=\phi(k)$ for $k=1,...,n$ and on 8  parameters of the monodromy matrix $M$,
while the $PSL(3,\R)$-equivalence  reduces the dimension by 8.
The pentagram map $T$ is generically defined  on the space ${\mathcal P}_n$. Namely,
for a twisted $n$-gon vertices of its image are the intersections of pairs of consecutive shortest diagonals:
$Tv_k:=(v_{k-1}, v_{k+1})\cap (v_{k}, v_{k+2})$. Such intersections are well defined for a generic point in  ${\mathcal P}_n$.
\medskip

a) {\it Results on integrability (in the twisted and closed cases).}
There is a Poisson structure on ${\mathcal P}_n$ invariant with respect to the pentagram map.
There are $2\lfloor n/2 \rfloor +2$ integrals in involution, which provide integrability of the pentagram map on  ${\mathcal P}_n$.
Its symplectic leaves have codimensions 2 or 4 in ${\mathcal P}_n$ depending on whether $n$ is odd or even,
and the invariant tori have dimensions $n-1$ or $n-2$, respectively \cite{OST99}.

Moreover, when restricted to the space ${\mathcal C}_n$ of closed polygons ($\dim {\mathcal C}_n=2n-8=2(n-4)$),
the map is still integrable and has invariant tori of dimension $n-4$ for odd $n$ and $n-5$ for even $n$.
Note that the space ${\mathcal C}_n$ of closed polygons is not a Poisson
submanifold in the space ${\mathcal P}_n$ of twisted $n$-gons,
so the corresponding Poisson structure on ${\mathcal P}_n$ cannot be restricted to ${\mathcal C}_n$.

There is a Lax representation for the pentagram map.
Coefficients of the corresponding spectral curve are the first integrals of the dynamics.
The pentagram map defines a discrete motion on the Jacobian of the spectral curve.
This motion is linear or staircase-like depending on the parity of $n$, see  \cite{FS11}.

\medskip

b) {\it Coordinates on  ${\mathcal P}_n$.}
The following two systems of coordinates on ${\mathcal P}_n$ are particularly convenient to work with, see \cite{OST99}.
Assume that $n$ is not divisible by 3.
Then there exists a unique lift of points $v_k=\phi(k)\in \RP^2$ to the vectors $V_k\in \R^3$ satisfying the condition
$\det|V_j, V_{j+1}, V_{j+2}|=1$ for each $j$. Associate a difference equation to a sequence of vectors $V_k\in \R^3$ by setting
$$
V_{j+3}=a_j V_{j+2} + b_j V_{j+1} + V_{j}
$$
for all $j\in \Z$.  The sequences $(a_j)$ and $(b_j)$ turn out to be $n$-periodic,
which is a manifestation of the fact that the lifts satisfy the relations $V_{j+n}=MV_j,\; j \in \Z,$ for a certain monodromy matrix $M\in SL(3,\R)$.
 The variables $a_j, b_j, \, 0\le j \le n-1$ are coordinates on the space ${\mathcal P}_n$.

There exists another coordinate system on the space ${\mathcal P}_n$, which is more geometric.
Recall that the cross-ratio of 4 points in $\PP^1$ is given by
$$
[t_1,t_2,t_3, t_4]=\frac{(t_1-t_2)(t_3-t_4)}{(t_1-t_3)(t_2-t_4)},
$$
where $t$ is any affine parameter. Now associate to each vertex $v_i$ the following two numbers, which are the cross-ratios of two 4-tuples of points lying on the lines $(v_{i-2},v_{i-1})$ and $(v_{i+1},v_{i+2})$ respectively:
$$
x_i=
[v_{i-2},\, v_{i-1},\, ((v_{i-2},  v_{i-1})\cap( v_{i}, v_{i+1})),\,  ((v_{i-2}, v_{i-1})\cap(v_{i+2}, v_{i+2}))]
$$
$$
y_i=
[((v_{i-2}, v_{i-1})\cap (v_{i+1},  v_{i+2})), (( v_{i-1}, v_{i}) \cap (v_{i+1}, v_{i+2})),\,
v_{i+1},\, v_{i+2}]
$$
In these coordinates the pentagram map has the form
$$
T^*x_i=x_i\frac{1-x_{i-1}y_{i-1}}{1-x_{i+1}y_{i+1}}\,\qquad
T^*y_i=y_{i+1}\frac{1-x_{i+2}y_{i+2}}{1-x_{i}y_{i}}\,.
$$
One can see that the pentagram map commutes with the scaling transformation \cite{OST99, Sch08}:
$$
R_s: \,(x_1,y_1, ..., x_n, y_n)\to (sx_1, s^{-1}y_1, ..., sx_n, s^{-1} y_n)\,.
$$
In these coordinates the invariant Poisson structure has a particularly simple form, see \cite{OST99}.

\medskip

c) {\it Continuous limit: the Boussinesq equation.}
The $n\to\infty$ continuous limit of a twisted $n$-gon with a fixed monodromy $M\in PSL(3,\R)$
can be viewed as a smooth parameterized curve $\gamma:\R\to \RP^2$  satisfying $\gamma(x+2\pi)=M\gamma(x)$ for all $x\in \R$.
The genericity assumption that every three consecutive points of an $n$-gon are in general position corresponds to the assumption
that $\gamma$ is a non-degenerate curve in $\RP^2$, i.e.,
the vectors $\gamma'(x)$ and $\gamma''(x)$ are linearly independent for all $x\in \R$.

\begin{figure}[hbtp]
\centering
\includegraphics[width=3in]{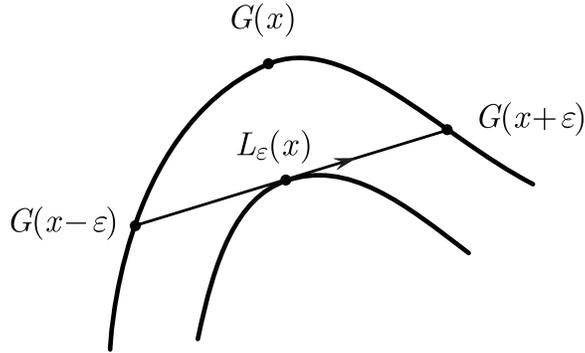}
\caption{\small Constructing the envelope $L_\ep(x)$ in 2D.}
\end{figure}

The space of such projectively equivalent  curves is
in one-to-one correspondence with linear differential operators
of the third order: $L=\partial^3+u_1(x)\partial+u_0(x)$, where the coefficients $u_0$ and $u_1$ are periodic in $x$.
Namely, a curve $\gamma(x)$ in $\RP^2$ can be lifted to a quasi-periodic curve $G=\{G(x)\}$ in
$\R^3$ satisfying $\det |G(x), G'(x), G''(x)|=1$ for all $x\in \R$. The components of the vector function $G(x)=(G_1(x),G_2(x),G_3(x))$
are homogenous coordinates of $\gamma(x)$ in $\RP^2$: $\gamma(x)=(G_1:G_2:G_3)(x)\in \RP^2$.
The vector function $G(x)$ can be identified with a solution of the unique linear differential operator $L$, i.e.,
the components of $G(x)$ are identified with three linearly independent solutions of the differential equation $Ly=0$.

\medskip

A continuous analog of the pentagram map is obtained by the following construction.
Given a non-degenerate curve $\gamma(x)$,  we draw the chord
$(\gamma(x-\ep),\gamma(x+\ep))$ at each point $x$. Consider the envelope $\ell_\ep(x)$ of these chords.
(Figure 2 shows their lifts: chords $(G(x-\ep),G(x+\ep))$ and their envelope $L_\ep(x)$.)
Let $u_{1,\ep}$ and $u_{0,\ep}$ be the  periodic coefficients of the corresponding differential operator. Their expansions in $\ep$ have the form
$u_{i,\ep}=u_i +\ep^2 w_i +{\mathcal O}(\ep^3)$ and allow one to define the evolution
$du_i/dt:=w_i$, $i=0,1$. After getting rid of $u_0$ this becomes the classical Boussinesq equation on the periodic function $u=u_1$, which is the $(2,3)$-flow in the KdV hierarchy of integrable equations on the circle: $u_{tt}+2(u^2)_{xx}+u_{xxxx}=0$.

Below we  generalize these results to higher dimensions.

\bigskip


\section{Geometric definition of higher pentagram maps}

\subsection{Pentagram map in  3D}

First we extend the notion of a closed polygon to a twisted one, similar to the 2D case.
We present the 3D case first, before giving the definition of the pentagram map in arbitrary dimension, since it is used in many formulas below.

\begin{definition}
{\rm
A {\it twisted $n$-gon in} $ \RP^3$ with a monodromy ${ M} \in SL(4,\R)$
is a map $\phi: \Z \to \RP^3$, such that
$\phi(k+n) =  M \circ \phi(k)$ for each $k\in \Z$.
(Here we consider the natural action of $SL(4,\R)$ on the corresponding projective space
$\RP^3$.)
Two twisted $n$-gons are (projectively) {\it equivalent} if there is a transformation $g \in SL(4,\R)$, such that $g \circ
\phi_1 = \phi_2$.
}
\end{definition}

Note that equivalent  $n$-gons must have similar monodromies.
Closed $n$-gons (space polygons) correspond to the monodromies $M=\mathrm{Id}$ and $-\mathrm {Id}$.
Let us assume that vertices of an  $n$-gon are in general position, i.e., no four consecutive vertices belong to one and the same plane in $\RP^3$.
Also, assume that $n$ is odd.
Then one can show (see Section \ref{sect:Neven} below and cf. Proposition 4.1 in \cite{OST99}) that there exists a unique lift of the vertices $v_k:=\phi(k) \in \RP^3$
to the vectors $V_k \in \R^4$ satisfying for all $j\in \Z$ the identities
$\det|V_j, V_{j+1}, V_{j+2}, V_{j+3}|=1$ and $V_{j+n}=MV_j$, where
$M\in SL(4,\R)$.\footnote{This explains our choice of the group $SL(4,\R)$ rather than the seemingly more natural group  $PSL(4,\R)$: since $SL(4,\R)$ is a two-fold cover of  $PSL(4,\R)$,
a twisted $n$-gon would have two different lifts from $\RP^3$ to $ \R^4$ corresponding to two different lifts of the monodromy $M$ from the latter group.}
These vectors satisfy difference equations
$$
V_{j+4} = a_j V_{j+3} + b_j V_{j+2} + c_j V_{j+1} - V_j, \; j \in \Z,
$$
with $n$-periodic coefficients $(a_j , b_j, c_j)$
and we employ this equation to introduce the $(a,b,c)$-coordinates on the space of twisted $n$-gons.

\medskip

Define the pentagram map on the classes of equivalent $n$-gons in such a way.

\begin{figure}[hbtp]
\centering
\includegraphics[width=2.9in]{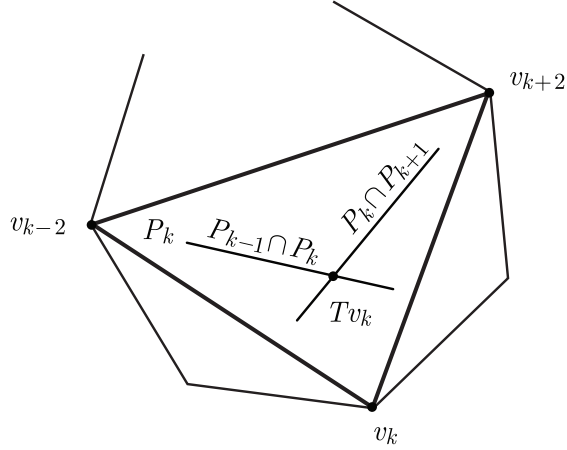}
\caption{\small In 3D the image $Tv_k$ of the vertex $v_k$  is the intersection of three ``short-diagonal" planes  $P_{k-1}, P_k, $ and $P_{k+1}$.}
\end{figure}

\begin{definition}
{\rm
Given an $n$-gon $\phi$ in  $\RP^3$, for each $k\in \Z$ consider the two-dimensional ``short-diagonal plane" $P_k:=(v_{k-2}, v_k, v_{k+2})$
passing through 3 vertices $v_{k-2}, v_k, v_{k+2}$. Take the intersection point of the three consecutive planes $P_{k-1}, P_k, P_{k+1}$
and call it the image of the vertex $v_k$ under the {\it space pentagram map} $T$, see Figure~3.
(We assume the general position, so that  every three consecutive planes $P_k$ for the given $n$-gon intersect at a point.)
}
\end{definition}

By lifting a two-dimensional plane $P_k$ from $\RP^3$ to the three-dimensional plane through the origin in $\R^4$ (and slightly abusing notation) we have $P_k =*(V_{k-2} \wedge V_k \wedge V_{k+2}) $
in terms of the natural duality $*$ between $\R^4$ and  $\R^{4*}$.
The lift of $Tv_k$ to $\R^4$ is proportional to $*[P_{k-1}  \wedge P_k  \wedge  P_{k+1}]$.

Below we describe the properties of this space pentagram map in detail.

\medskip


\subsection{Pentagram map in any dimension}\label{sec:highT}

Before defining the pentagram map in $\RP^d$, recall that
$SL(d+1,\R)$ is a two-fold cover of $PSL(d+1,\R)$  for odd $d$ and coincides with
the latter for even $d$.

\begin{definition}
{\rm
A {\it twisted $n$-gon in} $\RP^d$ with a monodromy $M \in SL(d+1,\R)$
is a map $\phi: \Z \to \RP^d$, such that
$\phi(k+n) =  M \circ \phi(k)$ for each $k\in \Z$, and where $ M$ acts naturally on $\RP^d$.
}
\end{definition}

We define the $SL(d+1,\R)$-equivalence of $n$-gons as above, and assume the vertices
$v_k:=\phi(k)$ to be  in general position, i.e., in particular, no $d+1$ consecutive vertices of an $n$-gon  belong to one and the same $(d-1)$-dimensional plane in $\RP^d$.

\begin{remark}\label{diff-eq}
{\rm
One can show that there exists a unique lift of the vertices $v_k=\phi(k) \in \RP^d$
to the vectors $V_k \in \R^{d+1}$ satisfying $\det|V_j, V_{j+1}, ..., V_{j+d}|=1$ and $ V_{j+n}=MV_j,\; j \in \Z,$ where
$M\in SL(d+1,\R)$, if and only if the condition $gcd(n,d+1)=1$ holds.
The corresponding difference equations have the form
\begin{equation}\label{eq:difference_anyD}
V_{j+d+1} = a_{j,d} V_{j+d} + a_{j,d-1} V_{j+d-1} +...+ a_{j,1} V_{j+1} +(-1)^{d} V_j,\quad j \in \Z,
\end{equation}
with $n$-periodic coefficients in the index $j$.
This allows one to introduce coordinates $\{ a_{j,k} ,\;0\le j\le n-1, \; 1\le k\le d \}$ on the space of twisted $n$-gons in $\RP^d$.
}
\end{remark}

\medskip

For a generic twisted $n$-gon in $\RP^d$ one can define the ``short-diagonal" $(d-1)$-dimensional  plane  $P_k$ passing
through $d$ vertices of the $n$-gon by taking every other vertex starting at the point $v_k$, i.e., through the vertices
$v_k, v_{k+2}, ..., v_{k+2d-2}$.
For calculations, however, it is convenient to have the set of vertices ``centered" at $v_k$, and
then the definition becomes slightly different in the odd and even dimensional cases.

Namely, for  odd dimension $d=2\varkappa+1$
we consider the {\it short-diagonal hyperplane} $P_k$ through the $d$ vertices
$$
P_k:=(v_{k-2\varkappa}, v_{k-2\varkappa+2},...,v_k, ..., v_{k+2\varkappa})
$$
(thus including the vertex $v_k$ itself), while for even dimension $d=2\varkappa$ we take $P_k$  passing through the  $d$ vertices
$$
P_k:=(v_{k-2\varkappa+1}, v_{k-2\varkappa+3}, ..., v_{k-1},v_{k+1}, ..., v_{k+2\varkappa-1})
$$
(thus excluding the vertex $v_k$).

\begin{definition}
{\rm
The {\it higher pentagram map} $T$ takes a vertex $v_k$ of a generic twisted
$n$-gon  in  $\RP^d$
to the intersection point of the $d$ consecutive short-diagonal planes $P_i$ around $v_k$.
Namely, for odd $d=2\varkappa+1$ one takes the intersection of the planes
$$
Tv_k:=P_{k-\varkappa}\cap P_{k-\varkappa+1}\cap ...\cap P_k\cap ...\cap P_{k+\varkappa}\,,
$$
while for even $d=2\varkappa$ one takes
the intersection of the planes
$$
Tv_k:=P_{k-\varkappa+1}\cap P_{k-\varkappa+2}\cap ...\cap P_{k}\cap ...\cap P_{k+\varkappa}\,.
$$
The corresponding map $T$ is well defined on the equivalence classes of  $n$-gons in $\RP^d$.
}
\end{definition}

As usual, we invoke the generality assumption  to guarantee that every $d$ consecutive hyperplanes $P_i$
 intersect at one point in $\RP^d$. It turns out that the pentagram map
defined this way has a special scaling invariance, which  allows one to prove its integrability:

\medskip

\noindent
{\bf Theorem C. (= Theorem \ref{thm:lax_anyD})}
{\it The scale-invariant higher pentagram map is a discrete completely integrable system
on equivalence classes  of $n$-gons in $\RP^d$. It has an explicit  Lax representation with a spectral parameter.}
\medskip

\begin{remark}
{\rm
The pentagram map defined this way in 1D is the identity map. 
In the 2D case this definition was given in \cite{Schwartz}  and its integrability 
for twisted polygons was proved in \cite{OST99}, while for closed ones
 it was proved in \cite{FS11, OST11}.
}
\end{remark}

\begin{remark}
{\rm
One can also give an ``asymmetric definition" for planes $P_k$, where more general sequences
of $d$ vertices $v_{k_j}$ are used, and then $Tv_k$ is defined as the intersection of $d$
consecutive planes $P_k$. It turns out, however, that exactly this ``uniform" definition of diagonal
planes $P_k$,
where $P_k$ passes through every other vertex, leads to integrability of the pentagram map.

One of possible  definitions of the pentagram map discussed in \cite{Beffa} coincides with ours in 3D.
In that definition one takes the intersection of a (possibly asymmetric, but containing the vertex $v_k$)
plane $P_k$  and the segment $[v_{k-1}, v_{k+1}]$: for our centered  choice of $P_k$ this segment belongs to both planes $P_{k-1}$ and $P_{k+1}$. The definitions become different in higher dimensions.
}
\end{remark}

\medskip


\subsection{General pentagram maps and duality}\label{generalmap}

We define more general pentagram maps $T_{p,r}$
depending on two integral parameters in arbitrary dimension $d$.
These parameters specify the diagonal planes and which of them to intersect.

\begin{definition}
{\rm
For a generic twisted $n$-gon in $\RP^d$ one can define a {\it $p$-diagonal hyperplane}  $P_k$
as the one passing through $d$ vertices of the $n$-gon by taking every $p$th vertex starting at the point $v_k$, i.e.,
$$
P_k:=(v_k, v_{k+p}, ..., v_{k+(d-1)p})\,.
$$
The image of the vertex $v_k$ under the {\it general pentagram map} $T_{p,r}$
is defined by intersecting every $r$th out of the $p$-diagonal hyperplanes  starting with $P_k$:
$$
T_{p,r}v_k:=P_{k}\cap P_{k+r}\cap ...\cap P_{k+(d-1)r}\,.
$$
The corresponding map $T_{p,r}$ is considered on the space of equivalence classes of  $n$-gons in
$\RP^d$.
}
\end{definition}

For the higher pentagram map $T$ discussed in Section \ref{sec:highT}
one has $p=2$, $r=1$, and the indices in the definition of $P_k$ are ``centered" at $v_k$.
In other words, $T=T_{2,1}\circ Sh$, where $Sh$ stands for some shift in the vertex index.
Below we denote by $Sh$ any shift in the index without specifying the shift parameter.
Note that $T_{p,p}=Sh$.

\begin{theorem}
There is a duality between the general pentagram maps $T_{p,r}$ and $T_{r,p}$:
$$
T_{p,r}=T^{-1}_{r,p}\circ Sh\,.
$$
\end{theorem}

For example, the map $T_{1,2}$ in 2D is defined by extending the sides of a polygon and intersecting them with the ``second neighbouring" sides.
This corresponds exactly to passing from $T(P)$ back to
$P$ in Figure~1, i.e., it is the inverse of $T$ modulo the numeration of vertices.

\proof
To prove this theorem we introduce the following ``duality maps,'' cf. \cite{OST99}.

\begin{definition}
{\rm
Given a generic sequence of points $\phi(j) \in \RP^d, \; j \in \Z,$  and a nonzero integer $p$
we define the {\it sequence}  $\alpha_p(\phi(j))\in (\RP^d)^*$ as the plane
$$
\alpha_p(\phi(j)):=(\phi(j), \phi(j+p),..., \phi(j+(d-1)p))\,.
$$
}\end{definition}

The following proposition is straightforward for our definition of the general pentagram map.

\begin{proposition}\label{prop:gen_pent}
For every nonzero $p$ the maps $\alpha_p$ are involutions modulo index shifts (i.e., $\alpha_p^2=Sh$),
they commute with the shifts (i.e., $\alpha_p\circ Sh= Sh\circ\alpha_p$), and
 the general pentagram map is the following composition: $T_{p,r}=\alpha_r\circ\alpha_p\circ Sh$.
\end{proposition}

To complete the proof of the theorem we note that
$$
T^{-1}_{r,p}=(\alpha_p\circ\alpha_r\circ Sh)^{-1}=Sh^{-1}\circ \alpha_r^{-1}\circ\alpha_p^{-1}
=Sh \circ\alpha_r \circ Sh\circ \alpha_p\circ Sh=\alpha_r \circ \alpha_p\circ Sh =T_{p,r}\circ Sh,
$$
since $\alpha_p^{-1}=\alpha_p\circ Sh$ from the Proposition above, while $Sh^{-1}=Sh$ and $Sh\circ Sh=Sh$.
\proofend

Below we construct a Lax form for the higher pentagram maps, i.e., for the maps $T_{2,1}$ (and hence for $T_{1,2}$ as well) for any $d$.
Integrability of the pentagram map on a special class of the so-called corrugated twisted  polygons in $\RP^d$ was proved in \cite{GSTV},
which should imply the integrability of the pentagram map $T_{p,1}$ in 2D.
Then the above duality would also give integrability of $T_{1,p}$ in $\RP^2$, defined as the intersection of a pair of polygon edges whose numbers differ by $p$.
(One should also mention that for $p$ and $r$ mutually prime with $n$ one can get rid of one of the parameters by appropriately  renumbering vertices, at least for the closed $n$-gon case, cf. \cite{OST11} in 2D.
This reduces the study to the  map $T_{p,1}$ for some values of $n$.)
Complete integrability of general pentagram maps for other pairs $(p,r)$ in $\RP^d$ is a wide open problem.

\begin{problem}
Which of the general pentagram maps $T_{p,r}$ in $\RP^d$ are completely integrable systems?
\end{problem}



\section{Continuous limit of the higher pentagram maps}

\subsection{Definition of the continuous limit}\label{sect:def_cont}

In this section we consider the continuous limit of polygons and the pentagram map on them.
In the limit  $n\to\infty$ a twisted $n$-gon becomes a smooth quasi-periodic curve $\gamma(x)$ 
in $\RP^d$.
Its lift $G(x)$ to $\R^{d+1}$ is defined by the conditions: $i)$ the components of the vector function 
$G(x):=(G_1(x),...,G_{d+1}(x))$ provide homogeneous coordinates for 
$\gamma(x)=(G_1:...:G_{d+1})(x)$ in $\RP^d$, $ii)$ $\det|G,G',...,G^{(d)}|(x)=1$ 
for all $x\in \R$, and $iii)$ $G(x+2\pi)=MG(x)$ for a given $M\in  SL(d+1,\R)$. Then $G$ satisfies a linear differential equation of order $d+1$:
\begin{equation}\label{eq:diff_anyD_onG}
G^{(d+1)}+u_{d-1}(x)G^{(d-1)}+...+u_1(x)G'+u_0(x)G=0
\end{equation}
with periodic coefficients $u_i(x)$. Here and below $'$ stands for $d/dx$.

Let us consider the case of odd $d=2\varkappa+1$. Fix a small $\epsilon>0$.
A continuous analog of the hyperplane $P_k$ is the hyperplane $P_\epsilon(x)$ passing through
$d$ points $\gamma(x-\varkappa\epsilon),..., \gamma(x),...,\gamma(x+\varkappa\epsilon)$ of the curve $\gamma$.\footnote{
For a complete analogy with the discrete case, one could take the points
$\gamma(x-2\varkappa\epsilon),\break \gamma(x-2(\varkappa-1)\epsilon),..., \gamma(x),...,\gamma(x+2\varkappa\epsilon)$.
However, one can absorb the  factor 2 by rescaling $\epsilon \to 2 \epsilon$.}
Note that as $\epsilon\to 0$ the hyperplanes $P_\epsilon(x)$ tend to the osculating hyperplane of the curve $\gamma$ spanned
by the vectors $\gamma'(x), \gamma''(x),...,\gamma^{(d-1)}(x)$ at the point $\gamma(x)$.

Let $\ell_\epsilon (x)$ be the envelope curve for the family of hyperplanes $P_\epsilon(x)$  for a fixed $\epsilon$.
The envelope condition means that for each $x$ the point $\ell_\epsilon (x)$ and the derivative vectors $\ell'_\epsilon (x),...,\ell^{(d-1)}_\epsilon (x)$
belong to the plane $P_\epsilon(x)$. This means that the lift of $\ell_\epsilon (x)$ to $L_\epsilon (x)$ in $\R^{d+1}$
satisfies the system of $d=2\varkappa+1$ equations (see Figure~4):
\begin{equation}\label{eq:Lfor_even}
\det | G(x-\varkappa\epsilon), G(x-(\varkappa-1)\ep ), ... , G(x),..., G(x+\varkappa\epsilon),  L^{(j)}_\epsilon(x) |=0,\quad j=0,...,d-1.
\end{equation}

\begin{figure}[hbtp]
\centering
\includegraphics[width=3.5in]{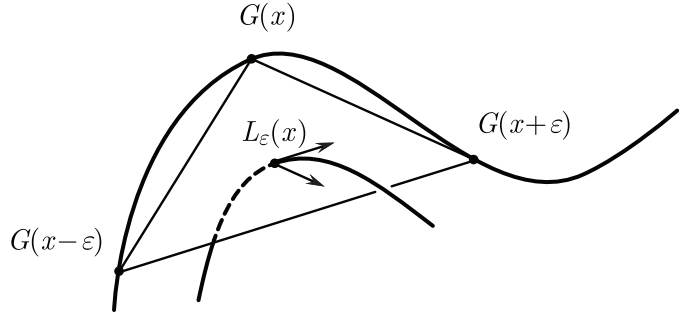}
\caption{\small The envelope $L_\ep(x)$ in 3D. The point $L_\ep(x)$ and the vectors $L_\ep'(x)$ and $L_\ep''(x)$ belong to the plane $(G(x-\ep), G(x), G(x+\ep))$.}
\end{figure}

Similarly, for even $d=2\varkappa$ the lift $L_\epsilon (x)$ satisfies the system of $d$ equations:
$$
\det | G(x-(2\varkappa-1)\epsilon), G(x-(2\varkappa-3)\ep ), ... , G(x-\ep),G(x+\ep),...
$$
\begin{equation}\label{eq:Lfor_odd}
 ...,
G(x+(2\varkappa-1)\epsilon),  L^{(j)}_\epsilon(x) |=0,\quad j=0,...,d-1.
\end{equation}

The evolution of the curve $\gamma$ in the direction of the envelope $\ell_\epsilon$, as $\epsilon$ changes,
defines a continuous limit of the pentagram map $T$. Namely, below we show that the expansion of $L_\epsilon(x)$ has the form
$$
L_\epsilon(x)=G(x)+\epsilon^2 B(x)+{\mathcal O} (\epsilon^4).
$$
The family of functions $L_\epsilon(x)$ satisfies a family of differential equations:
$$
L_\epsilon^{(d+1)}+u_{d-1,\ep}(x)L_\epsilon^{(d-1)}+...+u_{1,\ep}(x)L_\epsilon'+u_{0,\ep}(x)L_\epsilon=0, 
$$
where $u_{j,0}(x)=u_{j}(x).$

Expanding the coefficients $u_{j,\ep}(x)$ as $u_{j,\ep}(x)=u_{j}(x)+\ep^2w_j(x)+{\mathcal O}(\ep^4)$,
we define the continuous limit of the pentagram map $T$ as the system of evolution differential equations $du_j(x)/dt\, =w_j(x)$ for $j=0,...,d-1$,
i.e., $\ep^2$ plays the role of time.

\medskip

\noindent
{\bf Theorem A.}
{\it The continuous limit $du_j(x)/dt\, =w_j(x), \, j=0,...,d-1$ for $x\in S^1$ of the pentagram map is the $(2, d+1)$-KdV equation of the Adler-Gelfand-Dickey  hierarchy on the circle.}

\medskip

This theorem is proved as a combination of Theorems \ref{evolution-anyD} and \ref{thm:kdv} below.

\begin{remark}
{\rm
Recall the definition of the KdV hierarchy (after Adler-Gelfand-Dickey, \cite{Adler}).
Let $L$ be a linear differential operator of order $d+1$:
\begin{equation}\label{operatorL}
L = \partial^{d+1} + u_{d-1}(x) \partial^{d-1} + u_{d-2}(x) \partial^{d-2} + ...+ u_1(x) \partial + u_0(x)
\end{equation}
with periodic coefficients $u_j(x)$, where $\partial^{k}$ stands for $d^k/dx^k$.
One can define its fractional power
$L^{m/{d+1}}$ as a pseudo-differential operator for any positive integer $m$ and take
its purely differential part  $Q_m :=(L^{m/{d+1}})_+$ . In particular, for $m=2$ one has $Q_2= \partial^2 + \dfrac{2}{d+1}u_{d-1}(x) $. Then the $(m, d+1)$-KdV equation is the following evolution equation on (the coefficients of) $L$:
$$
\frac{d}{dt}{L} = [Q_m,L] \,.
$$
}
\end{remark}

\begin{remark}
{\rm
a) For $d=1$ the discrete pentagram map is the identity map, hence the continuous limit is trivial, which is consistent with vanishing of the (2,2)-KdV equation.

b) For $d=2$ the  (2,3)-KdV equation is the classical Boussinesq equation, found in \cite{OST99}.

c) Apparently, the $(2, d+1)$-KdV equation is a very robust continuous limit.
One obtains it not only for the pentagram map defined by taking every other vertex,
but also for a non-symmetric choice of vertices for the plane $P_k$, see Remark \ref{rm:diffchoice}.
Also, the  limit remains the same if instead of taking the envelopes one starts with a map defined by taking intersections of various planes  \cite{Beffa}. 
}
\end{remark}



\subsection{Envelopes and the KdV hierarchy}

\begin{theorem}\label{evolution-anyD}
For any dimension $d$, the envelope $L_\ep(x) $ has the expansion
$$
L_\ep(x) = G(x)+ {\ep^2} \,  C_d \left(  G''(x) +
\dfrac{2}{d+1}u_{d-1}G(x) \right) +{\mathcal O}(\ep^4)
$$
for a certain constant $C_d$, as $ \ep\to 0$.
\end{theorem}

The $\ep^2$-term of this expansion can be rewritten as $C_d\left(\partial^2+ \dfrac{2}{d+1}u_{d-1}(x)  \right)G(x).$
Consequently, it defines the following evolution of the curve $G(x)$:
$$
\frac{d}{dt} G= \left(\partial^2+ \dfrac{2}{d+1}u_{d-1}\right)G.
$$

\proof
Since $L_\epsilon$ approaches $G $ as $\epsilon\to 0$, one has the expansion
$L_\ep=G+\ep A+\ep^2 B+\ep^3 C+{\mathcal O}(\ep^4).$
First we note that  the expansion of $L_\ep$ in $\ep$ has only even powers of $\ep$, since
the equations (\ref{eq:Lfor_even}) and (\ref{eq:Lfor_odd}) defining $L_\ep$
have the symmetry $\epsilon\to - \ep$.
Therefore, we have $A=C=0$ and $L_\ep=G+\ep^2 B+{\mathcal O}(\ep^4).$

Notice that  $G(x)$ with its first $d$ derivatives form a basis in $\R^{d+1}$ for each $x$.
We express the vector coefficient $B$ in this basis: $B=b_0(x) G+b_1(x) G'+ ...+ b_d(x) G^{(d)}$.
Recall that, e.g., for odd  $d=2\varkappa+1$
the lift $L_\epsilon (x)$ satisfies the system of $d$ equations:
$$
\det |\, G(x-\varkappa\epsilon), G(x-(\varkappa-1)\ep ), ... , G(x),...,
G(x+\varkappa\ep),  L^{(j)}_\ep(x) \,|=0,\quad j=0,...,d-1.
$$
Fix $x$ and expand all terms in $\ep$: e.g., $G(x+\ep)=G(x)+\ep G'(x)+\frac{\ep^2}{2}G''+...$, etc.
In each equation consider the coefficients at the lowest power of $\ep$, being  $2+{d(d-1)}/{2}$ here.

The equation with $j=0$ gives $C_d\det |G, G',...,G^{(d-1)}, B|=0  $ for some nonzero $C_d$, which implies that there is no
$G^{(d)}$-term in the expansion of $B$. Similarly, for $j=1$ we obtain
$C_d\det |G, G',...,G^{(d-1)}, B'|=0  $, which means that  there is no
$G^{(d)}$-term in the expansion of $B'$, or, equivalently, there is no
$G^{(d-1)}$-term in the expansion of $B$.
Using this argument for $j=0,1,...,d-3,$ and $d-1$, we deduce that $B$ contains no
terms with $G^{(d)}, ..., G'''$ and $G'$.

The equation with $j=d-2$ results in a different term in the expansion and gives
$$
\det |G, G',...,G^{(d-1)}, B^{(d-2)}|=C_d\, \det |G, G',...,G^{(d-1)}, G^{(d)}| \,,
$$
which implies that $B=C_dG''+b(x)G$ for some function $b(x)$.

Finally, the normalization
$\det | \,L_\ep, L'_\ep, ..., L^{(d)}_\ep\,|=1$ allows one to find $b(x)$ by plugging in it
$L_\ep=G+\ep^2 (C_dG''+b(x)G)+{\mathcal O}(\ep^4).$ For the $\ep^2$-terms one obtains
$$
  C_d\det |\,G,G',...,G^{(d-2)},G^{(d+1)},G^{(d)}\,|
+ C_d\det |\,G,G',...,G^{(d-2)},G^{(d-1)},G^{(d+2)}\,|
$$
$$
+(d+1) b(x)\det |\,G,G',...,G^{(d-2)},G^{(d-1)},G^{(d)}\,|=0\,.
$$
By using the linear differential equation $G^{(d+1)}+u_{d-1}(x)G^{(d-1)}+...+u_0(x)G=0$
to express $G^{(d+1)}$ and $G^{(d+2)}$ via lower derivatives we see that the first
and the second determinants are equal to $-u_{d-1}$, while the last one is equal to 1.
Thus one has $(d+1) b(x)-2C_d\,u_{d-1}(x)=0$, which gives $b(x)=C_d\frac{2}{d+1}u_{d-1}(x)$.

Hence
$L_\ep=G+\ep^2 \,C_d\,(G''+\frac{2}{d+1}u_{d-1}G)+{\mathcal O}(\ep^4)$, as required.
\proofend

\begin{remark}\label{rm:diffchoice}
{\rm
One can see that the only condition on vertices defining the hyperplane $P_k$ required for the proof above
is that they are distinct.
A different choice of vertices for the hyperplane $P_k$ changes the constant $C_d$, but does not affect the evolution equation for $G$.
The above theorem  for an envelope  $L_\ep(x) $  is  similar to  an analogous expansion  in \cite{Beffa}  for a curve defined via certain plane intersections.
}
\end{remark}

\begin{theorem}\label{thm:kdv}
In any dimension $d$ the continuous limit of the pentagram map defined by the evolution
$$
\frac{d}{dt} G= \left(\partial^2+ \dfrac{2}{d+1}u_{d-1} \right)G
$$
of the curve $G$
 coincides with the $(2,d+1)$-KdV equation. Consequently, it is an infinite-dimensional completely integrable system.
\end{theorem}

\proof
Recall that the $(2,d+1)$-KdV equation is defined as the evolution
$$
\dot L=[Q_2, L]:=Q_2L-LQ_2,
$$
where the linear differential operator $L$ of order $d+1$ is given by formula (\ref{operatorL})
and $Q_2=\partial^2+ \dfrac{2}{d+1}u_{d-1} $. Here $\dot L$ stands for $dL/dt$.

By assumption, the evolution of the curve $G$  is described by the
differential equation $\dot G =Q_2G$. We would like to find
the evolution of the operator $L$ tracing it.
For any $t$, the curve $G$ and the operator $L$ are related by the differential equation $LG=0$ of order $d+1$.
Consequently,
$$
\dot L G +  L \dot G=0.
$$
Note that if the operator $L$ satisfies the $(2,d+1)$-KdV equation
and $G$ satisfies $\dot G =Q_2G$, we have the identity:
$$
\dot L G +  L \dot G=(Q_2L-LQ_2) G + L Q_2G= Q_2LG=0\,.
$$
In virtue of the uniqueness of the linear differential operator of the form  (\ref{operatorL})
for a given fundamental set of solutions, we obtain that indeed the evolution
of $L$ is described by the $(2,d+1)$-KdV equation.
\proofend

\begin{remark}
{\rm
The proof above is reminiscent of the one used in \cite{OK} to
study symplectic leaves of the Gelfand-Dickey brackets. Note that
the absence of the term linear in $\ep$ is related to the symmetric
choice of vertices for the hyperplane $P_k$. For a non-symmetric choice the
evolution would be defined by the linear term in $\ep$ and given by the
equation $\dot G=G'$. This is the initial, $(1,d+1)$-equation of the corresponding KdV hierarchy, manifesting the fact that the space $x$-variable  can be regarded
as the ``first time" variable.
A natural question arises whether the whole KdV-hierarchy is
hidden as an appropriate limit of the pentagram map.
An evidence to
this is given by noticing that the terms with higher powers in $\ep$ lead to
equations similar to the higher equations in the KdV hierarchy, see Appendix \ref{cont2D}.
}
\end{remark}

\begin{remark}
{\rm
One can see that the continuous limit of the general pentagram maps $T_{p,r}$  for various $p\not=r$ in $\RP^d$
defined via envelopes for a centered choice of vertices
is the same $(2,d+1)$-KdV flow, i.e., the limit depends only on the dimension.

Indeed, an analog of the $p$-diagonal is the plane $P_\ep(x)$
passing through the points \break
$G(x),  G(x+\ep p), ..., G(x+\ep(d-1)p)$. Rescaling $\ep$,
we can assume the points to be $G(x), G(x+\ep ), ..., G(x+\ep(d-1))$, which
leads to the planes  $P_\ep(x)$ defined in Section \ref{sect:def_cont} after a shift in $x$.
Then the definition of $L_\ep(x)$ via the envelope of such planes will give the same $(2,d+1)$-KdV equation.

It would be interesting to define the limit of the intersections of every $r$th plane
via some higher-order terms of the envelope, as mentioned in the above remark, so that it could lead to other $(m, d+1)$-equations in the KdV hierarchy.
}
\end{remark}
\bigskip



\section{Explicit formulas for the 3D pentagram map}

\subsection{Two involutions}\label{sect:invol}
Now we return to the 3D case.
In this section we assume that $n$ is odd and consider  $n$-gons in $\RP^3$.
Recall that in this case an $n$-gon with a given monodromy $M\in SL(4,\R)$
lifts uniquely to $\R^4$ and is described by difference equations
\begin{equation}\label{eq:abc_in3D}
V_{j+4} = a_j V_{j+3} + b_j V_{j+2} + c_j V_{j+1} - V_j,\quad j \in \Z,
\end{equation}
with $n$-periodic coefficients  $(a_j,b_j,c_j)$.
In other words, for odd $n$ the variables $(a_j,b_j,c_j), \; 0 \le j \le n-1,$ provide
coordinates on the space $\mathcal{P}_n$ of twisted $n$-gons in $\RP^3$ considered up to projective equivalence (see Proposition~\ref{abc-lemma}).

In order to find  explicit formulas for the pentagram map,
we present  it as a composition of two involutions
$\alpha$ and $\beta$, cf. Section \ref{generalmap},
and then find the formulas for each of them separately.
The same approach was used in \cite{OST99} in 2D, although the formulas in 3D are more complicated.

\begin{definition}
{\rm
Given a sequence of points $\phi(j) \in \RP^3, \; j \in \Z,$ define {\it two  sequences} $\alpha(\phi(j))\in (\PP^3)^*$
and $\beta(\phi(j))\in (\RP^3)^*$, where

$a)$ $\alpha(\phi(j))$ is the plane $(\phi(j-1), \phi(j), \phi(j+1))$;

$b)$ $\beta(\phi(j))$ is the plane $(\phi((j-2), \phi(j), \phi(j+2))$.
}
\end{definition}

\begin{proposition}
The maps $\alpha$ and $\beta$ are involutions, i.e., $\alpha^2=\beta^2=\rm{Id}$, while the pentagram map is their composition: $T=\alpha\circ\beta$.
\end{proposition}

Note that the indices which define the planes are symmetric with respect to $j$.
As a result, we do not have an extra shift of indices, cf. Proposition \ref{prop:gen_pent} (unlike the 2D case and the general map $\alpha_p$).

\begin{lemma}\label{map-alpha}
The involution $\alpha: V_i \to W_i = *(V_i \wedge V_{i-1} \wedge V_{i+1})$ maps
equation~(\ref{eq:abc_in3D}) to the following difference equation:
$$
W_{i+4} = c_{i+1} W_{i+3} + b_i W_{i+2} + a_{i-1} W_{i+1} - W_i.
$$
\end{lemma}

\begin{lemma}\label{map-beta}
The involution $\beta: V_i \to W_i =  \lambda_i *(V_i \wedge V_{i-2} \wedge V_{i+2})$ maps
equation~(\ref{eq:abc_in3D}) to the difference equation
$$
W_{i+4} = A_i W_{i+3} + B_i W_{i+2} + C_i W_{i+1} - W_i,
$$
where the coefficients $A_i,B_i,C_i$ are defined as follows:
\begin{align*}
A_i &= c_{i-1} (a_i a_{i+2} + a_{i+2} b_{i+1} c_i + c_i c_{i+2})^2 \lambda_{i+1} \lambda_{i+2} \lambda_{i+4}^2,\\
B_i &= ((a_{i-2}+b_{i-1}c_{i-2})(c_{i+2}+a_{i+2}b_{i+1})-a_{i+2}c_{i-2}) \times\\
    &\times (a_{i-1} a_{i+1} + a_{i+1} b_i c_{i-1} + c_{i-1} c_{i+1}) \lambda_i \lambda_{i+1} \lambda_{i+3} \lambda_{i+4},\\
C_i &= a_{i+1} (a_i a_{i+2} + a_{i+2} b_{i+1} c_i + c_i c_{i+2})^2 \lambda_{i+2} \lambda_{i+3} \lambda_{i+4}^2.
\end{align*}

The sequence $\lambda_i,\; i \in \Z,$ is $n$-periodic and is uniquely determined by the condition
$$
\lambda_i \lambda_{i+1} \lambda_{i+2} \lambda_{i+3} =
\dfrac{1}{(a_{i-2} a_i + a_i b_{i-1} c_{i-2} + c_{i-2} c_i)(a_{i-1} a_{i+1} + a_{i+1} b_i c_{i-1} + c_{i-1} c_{i+1})}.
$$
\end{lemma}

The proofs of these lemmas are straightforward computations, which we omit.
Combined together, these lemmas provide formulas for the pentagram map.
They have, however, one  drawback: one needs to solve a system of equations in $\lambda_i,\; i \in \Z,$
which results in the non-local character of the formulas in $(a,b,c)$-coordinates and their extreme complexity.

\medskip


\subsection{Cross-ratio type coordinates}

Similarly to the 2D case, there exist alternative coordinates on the space $\mathcal{P}_n$.
They are defined for any $n$, and the formulas for the pentagram map become local,
i.e., involving the vertex $\phi(j)$ itself and several neighboring ones.

\begin{definition}\label{def:x-y-z}
{\rm For odd $n$ the variables
$$
x_j = \dfrac{b_{j+1}}{a_j a_{j+1}}, \quad y_j =\dfrac{a_j}{b_{j+1} c_j}, \quad
z_j =\dfrac{c_{j+1}}{a_{j+1} b_j}, \quad 0 \le j \le n-1,
$$
provide coordinates on the space $\mathcal{P}_n$,
where the $n$-periodic variables $(a_j,b_j,c_j),\; j \in \Z,$ are defined by the difference equation
(\ref{eq:abc_in3D}).
}
\end{definition}

It turns out that the variables $(x_j,y_j,z_j),\; 0 \le j \le n-1$ are well defined and independent for
any $n$, even or odd. Below we provide two (equivalent) ways to define them for even $n$:
a pure geometric (local) definition of these variables (see Proposition \ref{prop:local-def})
and the above definition  extended to quasi-periodic sequences  $(a_j,b_j,c_j)$ (see Section \ref{sect:Neven}).





 \begin{theorem}\label{explicit}
 In the coordinates $x_i, y_i, z_i$  the pentagram map for any (either odd or even) $n$ is given by the formulas:
  \begin{align*}
 T^*(x_i) &= x_{i+1} \dfrac {1+y_{i-1}+z_{i+2}+y_{i-1}z_{i+2}-y_{i+1} z_i}{1+y_{i-1}+z_i},\\
 T^*(y_i) &= \dfrac{x_{i-1} y_{i-1}z_i}{x_i z_{i-1}} \dfrac{(1+y_{i+1}+z_{i+2})(1+y_{i-2}+z_{i-1})}{(1+y_i+z_{i+1})(1+y_{i-1}+z_{i+2}+y_{i-1} z_{i+2}-y_{i+1}z_i)},\\
 T^*(z_i) &= \dfrac{x_{i+1} z_i}{x_i}\dfrac{(1+y_{i+1}+z_{i+2})(1+y_{i-2}+z_{i-1})}{(1+y_{i-1}+z_i)(1+y_{i-2}+z_{i+1}-y_i z_{i-1}+y_{i-2}z_{i+1})}.
 \end{align*}
\end{theorem}

 Before proving this theorem  we describe the $(x_i,y_i,z_i)$ coordinates in greater detail.
It turns out that they may be defined completely independently of $(a_i,b_i,c_i)$ in the following geometric way.

\smallskip

Recall that the $x,y$ coordinates for the 2D pentagram map are defined as cross-ratios for quadruples of points on the line $(V_i, V_{i+1})$,
where two points are these vertices themselves, and two others are intersections of this line with extensions of the neighbouring edges.
Similarly, the next proposition describes the new coordinates via cross-ratios of
quadruples of points, $2$ of which are the vertices of an edge, and $2$ others are the intersection of the   edge  extension with two planes.
For instance, the variable
$y_i$ is the cross-ratio of 4 points on the line  $(V_i, V_{i+1})$, two of which are $V_i$ and $V_{i+1}$,
while two more are constructed as intersections of this line
with the planes via the triple $(V_{i+2}, V_{i+3}, V_{i+4})$ and with the plane via the triple $(V_{i+2}, V_{i+4}, V_{i+5})$. More precisely, the following proposition holds.

\begin{proposition}\label{prop:local-def}
The coordinates $x_i,y_i,z_i$ are given by the cross-ratios:
\begin{align*}
x_i &= -[V_{i+4},V_{i+5},\Phi^{45}_{012},\Phi^{45}_{123}],\\
y_i &= -[V_i,V_{i+1},\Phi^{01}_{234},\Phi^{01}_{245}],\\
z_i &= -[V_{i+4},V_{i+5},\Phi^{45}_{013},\Phi^{45}_{123}],
\end{align*}
where the point $\Phi^{j_1,j_2}_{m_1,m_2,m_3}$ for a given $i$ is the intersection of the line $(V_{i+j_1},V_{i+j_2})$
with the plane $(V_{i+m_1},V_{i+m_2},V_{i+m_3})$.
\end{proposition}

By the very definition these coordinates are projectively invariant.

\proofprop
If $*$ is the Hodge star operator with respect to the Euclidean metric in $\R^4\,(\,\mathrm{ or~ }\C^4)$,
then
$$
\Phi^{j_1,j_2}_{m_1,m_2,m_3}=
*\left( *(V_{i+j_1} \wedge V_{i+j_2}) \wedge *(V_{i+m_1} \wedge V_{i+m_2} \wedge V_{i+m_3}) \right).
$$
It suffices to prove the proposition in the case of an odd $n$, because the formulas are local,
and we can always add a vertex to change the parity of $n$.
Therefore, we may assume that $(a_j,b_j,c_j)$ are global coordinates and use them for the proof.

A simple computation shows that
\begin{align*}
\Phi^{45}_{012} &= -(b_{i+1} + a_i a_{i+1}) V_{i+4} + a_i V_{i+5},\\
\Phi^{01}_{234} &= V_i - c_i V_{i+1},\\
\Phi^{01}_{245} &= -b_{i+1} V_i + (a_i+c_i b_{i+1})V_{i+1},\\
\Phi^{45}_{013} &= (c_{i+1}+b_i a_{i+1})V_{i+4}-b_i V_{i+5},\\
\Phi^{45}_{123} &= -a_{i+1} V_{i+4} + V_{i+5}.
\end{align*}
Recall (see Lemma 4.5 in~\cite{OST99}) that if $4$ vectors $a,b,c,d \in \R^4\,(\mathrm{ or~ }\C^4)$ lie in the same 2-dimensional plane and are such that
$$
c = \lambda_1 a + \lambda_2 b, \quad d = \mu_1 a + \mu_2 b,
$$
then the cross-ratio of the lines spanned by these vectors in the plane is given by
$$
[a,b,c,d] = \dfrac{\lambda_2 \mu_1 - \lambda_1 \mu_2}{\lambda_2 \mu_1}.
$$
Comparing the cross-ratios with the original definition of the variables $x_i,y_i,z_i$ concludes the proof.
\proofend

Now we are in a position to prove the explicit local formulas in Theorem \ref{explicit}.
The proof is similar, but more involved than that of Proposition 4.11 in~\cite{OST99}.

\proofthm
Due to the local character of the  formulas for the pentagram map $T$ in $(x_j,y_j,z_j)$-coordinates, we may always add an extra vertex to make the number $n$ of vertices odd, and  then use coordinates $(a_j,b_j,c_j)$ and Lemmas~\ref{map-alpha} and~\ref{map-beta} for the proof.

The pentagram map is a composition $T = \alpha \circ \beta:  V_i \to U_i$. Namely,
$$
U_i = \mu_i *\left[ *(V_{i-3} \wedge V_{i-1} \wedge V_{i+1}) \wedge *(V_{i-2} \wedge V_i \wedge V_{i+2}) \wedge *(V_{i-1} \wedge V_{i+1} \wedge V_{i+3}) \right],
$$
where the constants $\mu_i$ are chosen so that $\det|U_j, U_{j+1}, U_{j+2}, U_{j+3}|=1$ for all $j$.

At the level of the coordinates $(a_j,b_j,c_j)$, we obtain:
$$
T^*(x_i) = \dfrac{B_{i+1}}{C_{i+1} C_{i+2}}, \quad
T^*(y_i) = \dfrac{C_{i+1}}{B_{i+1} A_{i-1}}, \quad
T^*(z_i) = \dfrac{A_i}{C_{i+2} B_i},
$$
where $A_i,B_i,C_i$ are defined in Lemma~\ref{map-beta}.
Eliminating the variables $\lambda_i$ with different $i$ by using the formula
for the product $\lambda_i \lambda_{i+1} \lambda_{i+2} \lambda_{i+3}$ concludes the proof.
\proofend

\medskip


\subsection{Coordinates on twisted polygons: odd vs. even $n$}\label{sect:Neven}

In this section we compare how one  introduces the coordinates on the space $\mathcal{P}_n$ of twisted $n$-gons for odd or even $n$,
and how this changes the statements above.

\begin{definition}\label{def:quasi-abc}
{\rm
Call a sequence $(a_j,b_j,c_j),\; j \in \Z$, $n$-{\it quasiperiodic} if
there is a sequence $t_j,\; j \in \Z$, satisfying
$t_j t_{j+1} t_{j+2} t_{j+3}=1$  and such that
\begin{equation}\label{abc-transform}
a_{j+n} = a_j \dfrac{t_j}{t_{j+3}},\quad b_{j+n} = b_j \dfrac{t_j}{t_{j+2}},\quad c_{j+n} = c_j \dfrac{t_j}{t_{j+1}}
\end{equation}
for each $j \in \Z$.
}
\end{definition}

Note that a sequence $t_j,\; j \in \Z,$ must be 4-periodic, and it is defined by  three parameters, e.g., by   $\alpha:=t_0/t_3, \,\beta:=t_0/t_2,$
and $ \gamma:=t_0/t_1$ with $\alpha\beta\gamma>0$, and hence $t_0=(\alpha\beta\gamma)^{1/4}$.
Thus the space ${\mathcal QS}_n$ of $n$-quasiperiodic sequences  has dimension
$3n+3$, and $\{(a_j,b_j,c_j),\; j =0,...,n-1\}\times (\alpha, \beta, \gamma)$ are coordinates on it.

Now we associate a sequence of vectors $V_j \in \R^4,\; j \in \Z,$ and
difference equations
\begin{equation}\label{diff-abc}
V_{j+4} = a_j V_{j+3} + b_j V_{j+2} + c_j V_{j+1} - V_j,\quad j \in \Z,
\end{equation}
to each twisted $n$-gon $v_j:=\phi(j) \in \RP^3, j \in \Z,$ with a monodromy $M\in SL(4,\R)$.
This gives a correspondence between sequences $(a_j,b_j,c_j), \; j \in \Z,$ and twisted $n$-gons.

\begin{proposition}\label{abc-lemma}
There is a one-to-one correspondence between twisted $n$-gons
(defined up to projective equivalence) and three-parameter equivalence classes in the space
${\mathcal QS}_n$ of $n$-quasiperiodic sequences $\{(a_j,b_j,c_j),\; j =0,...,n-1\}\times (\alpha, \beta, \gamma)$.

If $n$ is odd,  there exists a unique $n$-periodic sequence $(a_j,b_j,c_j)$ in each class.

If $n=4p$, then the numbers  $ \alpha, \beta , \gamma$ are projective invariants of a twisted $n$-gon.

If $n=4p+2$, then there is one projective invariant: $\alpha \gamma/\beta$.
\end{proposition}

In other words, for odd $n$ the equivalence classes are ``directed along" the parameters $(\alpha, \beta, \gamma)$
and one can chose a representative with $\alpha=\beta= \gamma=1$ in each class.
For $n=4p$ the classes are ``directed across" these parameters, and hence the latter are fixed for any given class.
The case $n=4p+2$ is in between: in a sense, two of the $(\alpha, \beta, \gamma)$-parameters and one of the $(a,b,c)$-coordinates
can serve as coordinates on each equivalence class.

This proposition can be regarded as an analogue of Proposition 4.1 and Remark 4.4 in \cite{OST99} for $d=2$.

\proof
First, we construct the correspondence, and then consider what happens for different arithmetics of $n$.
For a given $n$-gon $v_k=\phi(k)$ we construct a sequence of vertices $V_j \in \R^4,\; j \in \Z,$  in the following way: choose
the lifts $\phi(0) \to V_0, \;\phi(1) \to V_1,\; \phi(2) \to V_2$ arbitrarily, and then
determine the vectors $V_j,\; j > 2,$ and $V_j,\; j < 0,$ recursively using the condition $\det(V_j, V_{j+1}, V_{j+2}, V_{j+3})=1$, which follows from equation~(\ref{diff-abc}).

By definition of a twisted $n$-gon, we have $\phi(j+n) = M \circ \,\phi(j)$ for each $j \in \Z$,
where $M\in SL(4,\R)$. Consequently, for each $j \in \Z$ there exists a number $t_j$,
such that $V_{j+n} = t_j M V_j$, and the matrix $M$ is independent of $j$.
The equation $\det(V_{j+n}, V_{j+1+n}, V_{j+2+n}, V_{j+3+n})=1$ implies that
$t_j t_{j+1} t_{j+2} t_{j+3}=1$
and  $t_{j+4}=t_j$  for each $j \in \Z$. In other words, the whole sequence $\{t_j\}$ is determined by
$t_1, t_2, t_3,$ and then $t_0=1/t_1t_2t_3$.
Quasiperiodic conditions~(\ref{abc-transform})
follow from the comparison of the equation
$$
V_{j+4+n} = a_{j+n} V_{j+3+n} + b_{j+n} V_{j+2+n} + c_{j+n} V_{j+1+n} - V_{j+n}
$$
with equation~(\ref{diff-abc}).

Now we rescale the initial three vectors:
$V_0\mapsto k_0V_0, \; V_1\mapsto k_1V_1, \;V_2\mapsto k_2V_2$,
where $k_0 k_1 k_2\not=0$. A different  lift of the three initial vectors corresponds
to a different sequence $\tilde{V}_j = k_j V_j$, where the sequence $k_j, j \in \Z,$
must also be $4$-periodic and satisfy $k_0 k_1 k_2 k_3=1$.
This rescaling gives the action of $(\R^*)^3$ on the space ${\mathcal QS}_n$ of $n$-quasiperiodic sequences.
By construction, the corresponding orbits (or equivalence classes) of sequences are in a bijection with  twisted $n$-gons.
The group $(\R^*)^3$ acts as follows:
$$
t_j \mapsto  t_j (k_j/k_{j+n}),\;
a_j\mapsto a_j (k_{j}/k_{j+3}), \;
b_j\mapsto b_j (k_{j}/k_{j+2}), \;
c_j\mapsto c_j (k_{j}/k_{j+1}).
$$

\smallskip

Now we have 3 cases:
\begin{itemize}
\item $n$ is odd.
Then the above $(\R^*)^3$-action on $t_0,\,t_1,\,t_2$ allows one to make them all equal to 1, which corresponds to the constant sequence $\{t_j=1,\,j\in\Z\}$ and
a periodic sequence $\{(a_j,b_j,c_j),\; j \in \Z\}$. Indeed, e.g., for $n=4p+3$ one has
the system of 3 equations: $t_0(k_0/k_3)=1, \,t_1(k_1/k_4)=1,\,t_2(k_2/k_5)=1$.
Since $k_4=k_0$, $k_5=k_1$, and $k_3=1/k_0k_1k_2$, we obtain a system of 3 equations
on the unknowns $k_0,k_1,k_2$, which has the unique solution.

\item $n=4p+2$.
One can check that the $(\R^*)^3$-action does not change the ratio $\dfrac{t_0 t_{2}}{t_{1} t_{3}}=\dfrac{\alpha \gamma}{\beta}$.

\item $n=4p$.
The $(\R^*)^3$-action does not change the three quantities $t_0/t_3, \, t_0/t_2, \, t_0/t_1$.
\proofend
\end{itemize}

\medskip

Now we can introduce coordinates on the space
${\mathcal QS}_n$ of $n$-quasiperiodic sequences using Definition \ref{def:x-y-z} and
quasiperiodic variables $(a_j,b_j,c_j)$.

\begin{proposition}\label{xyz-relations}
For any $n$ the variables $(x_j,y_j,z_j),\; 0 \le j \le n-1$ are independent and
constant on the equivalence classes in ${\mathcal QS}_n$, i.e., they are well-defined local
coordinates on the space $\mathcal{P}_n= {\mathcal QS}_n/\sim$.
\end{proposition}

\proof
It is straightforward to check that the $(\R^*)^3$-action defined above is trivial on the variables $(x_j,y_j,z_j)$.
For instance,
$$
 x_j =\dfrac{b_{j+1}}{a_j a_{j+1}} \to  \dfrac{b_{j+1}(k_{j+1}/k_{j+3})}{a_j (k_{j}/k_{j+3})\,a_{j+1}
(k_{j+1}/k_{j+4})}=x_j\,.
$$
The independence of the new variables on ${\mathcal QS}_n$ follows from that for the original ones.
Alternatively, it also follows from their local geometric definition (Proposition \ref{prop:local-def}).
\proofend

\begin{remark}
{\rm
In the $(a,b,c)$-coordinates for even $n$ some of the $(\alpha, \beta , \gamma)$-parameters were needed  to describe the equivalence classes in ${\mathcal QS}_n$. In the $(x,y,z)$-coordinates those parameters
are functions of $x_j,y_j,z_j$:

$i)$ for $n=4p+2$,
$$
\prod_{j=0}^{2p} \dfrac{x_{2j}^2 y_{2j} z_{2j+1}}{ x_{2j+1}^2 y_{2j+1} z_{2j}}
=\dfrac{\alpha \gamma}{\beta}\;;
$$

$ii)$ for $n=4p$,
$$
\prod_{j=0}^{p-1} \dfrac{x_{4j} x_{4j+2} y_{4j+2} z_{4j+1}}{x_{4j+1} x_{4j+3} y_{4j+3} z_{4j+2}}
= \alpha,\quad
\prod_{j=0}^{p-1} \dfrac{y_{4j+1} z_{4j}}{y_{4j+3} z_{4j+2}} =\beta,\quad
 \prod_{j=0}^{p-1} \dfrac{y_{4j} z_{4j+3}}{y_{4j+2} z_{4j+1}} =\frac \gamma\alpha\,.
$$
These identities follow from Definitions \ref{def:x-y-z} and \ref{def:quasi-abc}.}
\end{remark}


\section{Algebraic-geometric integrability of the 3D pentagram map}\label{sect:alg-geo}

In this section we complexify the pentagram map and assume that everything is defined over $\C$.

\subsection{Scaling transformations and a Lax function in 3D}\label{S:lax-func}

Recall that a discrete Lax equation with a spectral
 parameter is a representation of a dynamical system in the form
\begin{equation}\label{lax-eq}
 L_{i,t+1}(\lambda) = P_{i+1,t}(\lambda) L_{i,t}(\lambda) P_{i,t}^{-1}(\lambda),
 \end{equation}
where $t$ stands for the discrete time variable,  $i$ refers to the vertex index,
and $\lambda$ is a complex spectral parameter.

The pivotal property responsible for algebraic-geometric integrability of all pentagram maps
considered in this paper is the presence of a scaling invariance.
In the 2D case, this means the invariance with respect to the transformations
$a_j \to a_j s, \; b_j \to b_j/s$, where $s$ is an arbitrary number.
In the 3D case, the pentagram map is invariant with respect to the transformations
$a_j \to a_j s, \; b_j \to b_j, \; c_j \to c_j s$. In both cases the invariance follows from the explicit formulas of the map.
Note that formally one can define other pentagram maps by choosing the intersection planes in many different
ways. However, only very few of these maps  possess any scaling invariance.
Below we derive a Lax representation from the scaling invariance.
First we do it for odd $n$, when $(a_j,b_j,c_j),0 \le j \le n-1,$ are coordinates on the space ${\mathcal P}_n$.

\begin{theorem}\label{lax-abc}
The 3D pentagram map on twisted $n$-gons with  odd $n$ admits a Lax representation with the Lax function $L_j(\lambda)$ given by
 \[
 L_j(\lambda) =
 \begin{pmatrix}
  c_j/\lambda & 1/\lambda & 0 & 0\\
  b_j         & 0        & 1  & 0\\
  a_j/\lambda & 0        & 0  & 1/\lambda\\
  -1          & 0        & 0  & 0
 \end{pmatrix} =
 \begin{pmatrix}
  0       & 0   & 0       & -1\\
  \lambda & 0   & 0       & c_j\\
  0       & 1   & 0       & b_j\\
  0       & 0   & \lambda & a_j
 \end{pmatrix}^{-1}
 \]
  in the coordinates $a_j,b_j,c_j$. Its determinant is $\det{L_j} \equiv 1/\lambda^2$.
\end{theorem}

Note that we always consider a polygon and the corresponding Lax function at a particular
moment of time. Whenever necessary we  indicate the moment of time explicitly by adding the second index ``$t$''
to the Lax function (above $L_j:=L_{j,t}$), while if there is no ambiguity we keep only one index.
Before proving this  theorem we give the following

\begin{definition}
 {\rm
 {\it The monodromy operators} $T_{0,t},T_{1,t},...,T_{n-1,t}$ are defined as
 the following ordered products of the corresponding Lax functions:
 \begin{align*}
 &T_{0,t} = L_{n-1,t} L_{n-2,t} ... L_{0,t},\\
 &T_{1,t} = L_{0,t} L_{n-1,t} L_{n-2,t} ... L_{1,t},\\
&  ...\\
  &T_{i,t} = L_{i+n-1,t} L_{i+n-2,t} ... L_{i+1,t} L_{i,t},\\
 &  ...\\
 &T_{n-1,t} = L_{n-2,t} L_{n-3,t} ... L_{0,t} L_{n-1,t},
 \end{align*}
where the (integer) index $t$  represents the moment of time.
}
 \end{definition}

\proofthm
First observe that the Lax equation implies that the corresponding monodromy operators satisfy
$$
T_{i,t+1}(\lambda) = P_{i,t}(\lambda) T_{i,t}(\lambda) P_{i,t}^{-1}(\lambda),
$$
i.e., $T_{i,t}(\lambda)$ changes to a similar matrix when $t\to t+1$, and hence
 the eigenvalues of the matrices $T_{i,t}(\lambda)$ as functions of $\lambda$
are invariants of the map. Conversely, if some function $T_{i,t}(\lambda)$ has this property,
then there must exist a matrix $P_{i,t}(\lambda)$
(defined up to a multiplication by a scalar function) satisfying the above equation.

How to define such a monodromy depending on a parameter?
The monodromy matrix associated with the difference equation
$$
V_{j+4} = a_j V_{j+3} + b_j V_{j+2} + c_j V_{j+1} - V_j
$$
is $M = N_0 N_1 N_2... N_{n-1}$, where
\[
N_j =
 \begin{pmatrix}
  0 & 0 & 0 & -1\\
  1 & 0 & 0 & c_j\\
  0 & 1 & 0 & b_j\\
  0 & 0 & 1 & a_j
 \end{pmatrix}.
 \]
For odd $n$ the variables $(a_j,b_j,c_j), 0 \le j \le n-1,$ are well-defined coordinates on the space of twisted
$n$-gons. These variables are periodic: for any $j$ we have $a_{j+n}=a_j,\; b_{j+n}=b_j,\; c_{j+n}=c_j$.
The vectors $V_j$ are quasi-periodic: $V_{j+n}=M V_j$, and depend on the lift of the points from
the projective space. This means that the pentagram map preserves the eigenvalues of the matrix $M$,
but not the matrix $M$ itself.

Lemmas~\ref{map-alpha} and~\ref{map-beta} imply that the pentagram map is invariant
with respect to the scaling transformations: $a_j \to s a_j, \; c_j \to s c_j$.
Therefore, the pentagram map also preserves the eigenvalues of the monodromy
matrix $M(s)$ corresponding to the $n$-gons scaled by $s$.
Namely, we have
\[
M(s) = N_0(s) N_1(s) N_2(s)... N_{n-1}(s), \text{ where }
N_j(s) =
 \begin{pmatrix}
  0 & 0 & 0 & -1\\
  1 & 0 & 0 & s c_j\\
  0 & 1 & 0 & b_j\\
  0 & 0 & 1 & s a_j
 \end{pmatrix}.
\]
Now one can see that the matrix $N_j(s)$ can be chosen as a Lax function.
For technical reasons (which will be clear later), we define the Lax matrix as $L_j^{-1}(\lambda) := (g^{-1} N_j(s) g)/s$,
where $g := \text{diag}(1,s,1,s)$, and $\lambda := 1/s^2$. This gives the required matrix $L_j(\lambda)$.\proofend

\medskip

As we mentioned before, the formulas for the pentagram map are non-local in the variables
$(a_j,b_j,c_j)$. As a result, an explicit expression for the matrix $P_{i,t}(\lambda)$ becomes non-local as well.
On the other hand, one can use the variables $(x_j,y_j,z_j)$ (given by Definition \ref{def:x-y-z}) to describe a Lax representation.
Their advantage is that all formulas become local and are valid for any $n$, both even and odd.

\begin{theorem}\label{xyz-lax}
For any $n$ the equations for the 3D pentagram map  are equivalent to the Lax equation
$$
\tilde{L}_{i,t+1}(\lambda) = \tilde{P}_{i+1,t}(\lambda) \tilde{L}_{i,t}(\lambda) \tilde{P}_{i,t}^{-1}(\lambda),
$$
where
\[
\tilde{L}_{i,t}(\lambda) =
\begin{pmatrix}
  0               & 0   & 0           & -1\\
  \lambda x_i y_i & 0   & 0           & 1\\
  0               & z_i & 0           & 1\\
  0               & 0   & \lambda x_i & 1
 \end{pmatrix}^{-1},
\]
\[
\tilde{P}_{i,t}(\lambda)=
\begin{pmatrix}
  0 & \rho_i & 0 & -\rho_i\\
  \lambda \sigma_i(1+z_i) & -\rho_i & \lambda \sigma_i & \rho_i\\
  y_{i-1} \theta_i & \dfrac{z_{i-1}}{\tau_i} & -\theta_i & \dfrac{1+y_{i-2}}{\tau_i}\\
  -\dfrac{\lambda y_{i-1}}{1+y_{i-1}+z_i} & 0 & \dfrac{\lambda}{1+y_{i-1}+z_i} & 0
 \end{pmatrix},
\]
and the variables  $\rho_i,  \sigma_i, \tau_i,$ and $\theta_i$ stand for
\begin{align*}
\rho_i &= \dfrac{1}{x_i (1+y_i+z_{i+1})},\\
\sigma_i &= \dfrac{x_{i-1} y_{i-1}(1+y_{i-2}+z_{i-1})}{x_i z_{i-1}(1+y_{i-1}+z_i)(1+y_i+z_{i+1})},\\
\tau_i &= x_i(1+y_{i-2}-y_i z_{i-1}+z_{i+1}+z_{i+1} y_{i-2}),\\
\theta_i &=\dfrac{1+y_{i-2}+z_{i-1}}{\tau_i (1+y_{i-1}+z_i)}.
\end{align*}
\end{theorem}
\proof
The proof is a long but straightforward verification.
\proofend
\begin{remark}\label{lax-abc-xyz}
{\rm
The Lax functions $L$ and $\tilde{L}$ in the $(a_i,b_i,c_i)$ and $(x_i,y_i,z_i)$ variables are related to each other as follows:
$$
\tilde{L}_{i,t} = a_{i+1} \left( h_{i+1}^{-1} L_{i,t} h_i\right), \text{ where } h_i := \text{diag}(1,c_i,b_i,a_i).
$$
}\end{remark}

\medskip


\subsection{Properties of the spectral curve}

Recall that the monodromy operators $T_{i,t}(\lambda)$ satisfy the equation
$$
T_{i,t+1}(\lambda) = P_{i,t}(\lambda) T_{i,t}(\lambda) P_{i,t}^{-1}(\lambda).
$$
It implies that the function of two variables $R(\lambda,k)=\det{(T_{i,t}(\lambda) - k \,\text{Id})}$ is
independent of $i$ and $t$. Furthermore, $R(\lambda,k)=0$ is  a polynomial relation between $\lambda$ and $k$:  $R(\lambda,k)$ becomes a polynomial after a multiplication by a power of $\lambda$.
Its coefficients are integrals of motion for the pentagram map. The zero set of $R(\lambda,k)=0$
is an algebraic curve in $\C^2$. A standard procedure (of adding the infinite points and normalization with
a few blow-ups) makes it into a compact Riemann surface, which we call the {\it spectral curve} and denote by $\Gamma$.
In this section we explore some of the properties of the spectral curve and, in particular, find its genus.

\begin{definition}\label{ijg-def}
{\rm
For an odd $n$ define the {\it spectral function} $R(\lambda,k)$ as
$$
R(\lambda,k):=\det{(T_{i,t}(\lambda) - k \,\text{Id})},
$$
i.e., using the Lax function in the $(a,b,c)$-coordinates from Theorem~\ref{lax-abc}.
The {\it spectral curve} $\Gamma$ is the normalization of the compactification of the curve $R(\lambda,k)=0$.

We define the {\it integrals of motion} $I_j,J_j,G_j,\; 0 \le j \le q=\lfloor n/2 \rfloor,$ as
the coefficients of the expansion
$$
R(\lambda,k) = k^4 - k^3 \left( \sum_{j=0}^q G_j \lambda^{j-n} \right) + k^2 \left( \sum_{j=0}^q J_j \lambda^{j-q-n} \right)
- k \left( \sum_{j=0}^q I_j \lambda^{j-2n} \right) + \lambda^{-2n} = 0.
$$
}
\end{definition}

When $n$ is even, the sequence $(a_j,b_j,c_j),\; j \in \Z$, is not $n$-periodic,
and the monodromy operator $T_{i,t}(\lambda)$ cannot be defined. One should use the Lax function $\tilde{L}_{i,t}(\lambda)$ in the $(x,y,z)$-coordinates from
Theorem~\ref{xyz-lax} to define the monodromy operator $\tilde{T}_{i,t}(\lambda)$ and the corresponding spectral curve.

Namely, first note that the integral of motion  $I_0$ has the following explicit expression:
$$
I_0 = \prod_{i=0}^{n-1} a_i = \left( \prod_{i=0}^{n-1} x_i^2 y_i z_i \right)^{-1/4}.
$$

\begin{definition}\label{ijg-def-2}
{\rm
For any $n$ (either even or odd), the {\it spectral function} is
$$
R(\lambda,k) = \tilde{R}(\lambda, k I_0)/I_0^4,\; \text{ where }\;
 \tilde{R}(\lambda,k):=\det{(\tilde{T}_{i,t}(\lambda) - k \,\text{Id})}\,,
 $$
 and the monodromy operator $\tilde{T}_{i,t}(\lambda)$ is defined using the Lax function $\tilde{L}_{i,t}(\lambda)$
from Theorem~\ref{xyz-lax}.
}
\end{definition}

The spectral function $R(\lambda,k)$ defined this way coincides with $\det{(T_{i,t}(\lambda) - k \,\text{Id})}$ for odd $n$, since $\tilde{T}_{i,t} = I_0 \left( h_i^{-1} T_{i,t} h_i \right)$, see Remark~\ref{lax-abc-xyz}.
It is convenient to have such a unified definition for computations of  integrals of motion.
\medskip

\begin{theorem}\label{thm:genus}
For generic values of the integrals of motion $I_j,J_j,G_j$,
the genus $g$ of the spectral curve $\Gamma$ is $g = 3q$ for odd $n$ and $g = 3q-3$ for even $n$, where $q=\lfloor n/2 \rfloor$.
\end{theorem}

To prove it, we first describe the singularities of $R(\lambda,k)=0$ by considering the formal series solutions
(the so-called Puiseux series).

\begin{lemma}\label{spec-sing}
If $n$ is even, the equation $R(\lambda,k)=0$ has 4 formal series solutions at $\lambda=0$:
\begin{align*}
 O_1: \quad k_1 &= \dfrac{1}{I_0} - \dfrac{I_{1}}{I_0^2} \lambda + {\mathcal O}(\lambda^2),\\
 O_{2,3}: \quad k_{2,3} &= \dfrac{k_*}{\lambda^q} + {\mathcal O}\left(\dfrac{1}{\lambda^{q-1}}\right), \quad \text{ where $k_*$ satisfies} \quad G_0 k_*^2 - J_0 k_* + I_0 = 0,\\
 O_4: \quad k_4 &= \dfrac{G_0}{\lambda^n} + \dfrac{G_1}{\lambda^{n-1}}+\dfrac{G_2}{\lambda^{n-2}}+{\mathcal O}(\lambda^{3-n}),
\end{align*}
and $4$ solutions at $\lambda=\infty$:
$$
 W_{1,2,3,4}: \quad k_{1,2,3,4}= \dfrac{k_\infty}{\lambda^q} + {\mathcal O}\left( \dfrac{1}{\lambda^{q+1}} \right), \text{ where }
 k_\infty^4 - G_q k_\infty^3 + J_q k_\infty^2 - I_q k_\infty + 1 = 0.
$$
If $n$ is odd, the equation $R(\lambda,k)=0$ has 4 formal series solutions at $\lambda=0$:
\begin{align*}
 O_1: \quad k_1 &= \dfrac{1}{I_0} - \dfrac{I_1}{I_0^2} \lambda + {\mathcal O}(\lambda^2),\\
 O_2: \quad k_{2,3} &= \pm \dfrac{\sqrt{-I_0/G_0}}{\lambda^{n/2}} + \dfrac{J_0}{2 G_0 \lambda^{(n-1)/2}} + {\mathcal O}\left(\dfrac{1}{\lambda^{(n-2)/2}}\right),\\
 O_3: \quad k_4 &= \dfrac{G_0}{\lambda^n} + \dfrac{G_1}{\lambda^{n-1}}+\dfrac{G_2}{\lambda^{n-2}}+{\mathcal O}(\lambda^{3-n}),
\end{align*}
and $4$ solutions at $\lambda=\infty$:
$$
 W_{1,2}: \quad k_{1,2,3,4}= \dfrac{k_\infty}{\lambda^{n/2}} + {\mathcal O}\left( \dfrac{1}{\lambda^{(n+1)/2}} \right), \text{ where }
 k_\infty^4 + J_q k_\infty^2 + 1 = 0.
$$
The remaining coefficients of the series are determined uniquely.
\end{lemma}
\prooflem
One finds the series coefficients  recursively by substituting the series into the equation $R(\lambda,k)=0$, which determines the spectral curve.
\proofend

Now we can complete the proof of Theorem~\ref{thm:genus}.

\proofthm
As follows from the definition of the spectral curve $\Gamma$, it is a ramified 4-fold cover of $\CP^1$, since the $4\times 4$-matrix $\tilde{T}_{i,t}$ (or ${T}_{i,t}$) has 4 eigenvalues.
By the Riemann-Hurwitz formula the  Euler characteristic of $\Gamma$ is
$\chi(\Gamma)=4\chi(\CP^1)-\nu=8-\nu$, where $\nu $ is the number of  branch points.
On the other hand, $\chi(\Gamma)=2-2g$, and once we know $\nu$ it allows us to find the genus
of the spectral curve $\Gamma$ from the formula $2-2g=8-\nu$.

The number $\nu$ of branch points of $\Gamma$ on the $\lambda$-plane equals
the number of zeroes of the function $\partial_k R(\lambda,k)$
aside from the singular points.
The function $\partial_k R(\lambda,k)$ is meromorphic on $\Gamma$, therefore the number of its zeroes equals
the number of its poles.
One can see that for any $n$ the function $\partial_k R(\lambda,k)$ has poles of total order $9n$ at $z=0$,
and it has zeroes of total order $6n$ at $z=\infty$.
Indeed, substitute the local  series from Lemma \ref{spec-sing} to the expression for
$\partial_k R(\lambda,k)$.  (E.g., for $n=2q$ at $O_4$ one has $k\sim \lambda^{-n}$.
The leading terms of $\partial_k R(\lambda,k)$ for the pole at $\lambda=0$ are
$k^3, k^2\lambda^{-n}, k\lambda^{-q-n}, \lambda^{-2n}$. The  first two terms, being  of order $\lambda^{-3n}$, dominate and give the pole order of $3n=6q$.)
The corresponding orders of the poles and zeroes of $\partial_k R(\lambda,k)$ on  $\Gamma$ are summarized as follows:
\[
\begin{array}{||c|c||c|c||c|c||c|c||}
\hline
n=2q & \text{ pole } &     & \text{ zero } & n=2q+1 & \text{ pole } &     & \text{ zero }\\ \hline
O_1  &  4q           & W_1 & 3q            & O_1   &  2n           & W_1 & 3n\\ \hline
O_2  &  4q           & W_2 & 3q            & O_2   &  4n           & W_2 & 3n\\ \hline
O_3  &  4q           & W_3 & 3q            & O_3   &  3n           & & \\ \hline
O_4  &  6q           & W_4 & 3q            & & & &\\ \hline

\end{array}
\]

 For instance, for $n=2q$ this gives the total order of poles: $4q+4q+4q+6q=18q=9n$, while the total
order of zeroes is $4\times 3q=12q=6n$.

Therefore, the number of zeroes of $\partial_k R(\lambda,k)$ at nonsingular points $\lambda\not=\{0,\infty\}$  is $\nu=9n-6n=3n$, and so is
the total number of branch points of $\Gamma$ in the finite part of the $(\lambda,k)$ plane.
The difference between odd and even values of $n$ arises because $\Gamma$ has 2 additional branch points at $\lambda=\infty$,
and 1 branch point at $\lambda=0$ for odd $n$, i.e., $\nu=3n+3$.

Consequently,  one has $2-2g=8-\nu $ with $\nu=3n$ for even $n$ and $\nu=3n+3$  for odd $n$.
The required expression for the genus $g$ follows:
$g=3q-3$ for $n=2q$ and $g=3q$ for $n=2q+1$.
\proofend

\medskip

\begin{remark}
{\rm
Now we describe a few integrals of motion using the coordinates $(a_j,b_j,c_j), 0 \le j \le n-1,$ when $n$ is odd.
The description  is similar to that in the 2D case (cf. Section 5.2 and Proposition 5.3 in~\cite{OST99}).
Consider a code which is an ordered sequence of digits from $1$ to $4$. The number of digits in a code is
$p,q,r,t$, respectively. The code is called ``admissible'' if $p+2q+3r+4t=n$.
The number $p+r$ is called its ``weight.''
Each code expands in a ``word'' of $n$ characters in the following way: $1,2,3,4$ are replaced by
``a'',``*b'',``**c'',``****'', respectively. Now we label the vertices of an $n$-gon by $0,1,...,n-1$,
and associate letters in a word to them keeping the order. We obtain one monomial by taking the product of the variables
$a_i$,$b_i$, or $c_i$ that occur at the vertex $i$. The letter ``*'' corresponds to ``1''. The sign of the monomial
is $(-1)^t$. Next step is to permute the numbering of the vertices cyclically and take the sum of the monomials.
Note, however, that if, for example, $n=9$, then the code ``333'' corresponds to the sum
$c_0 c_3 c_6+c_1 c_4 c_7 + c_2 c_5 c_8$ without the coefficient $3$.
Finally, we repeat this procedure for all admissible codes of weight $p+r$ and denote the total sum by $\hat{I}_{p+r}$.
Additionally, we define the sum $\hat{G}_{p+r}$ by substituting $a_i \to c_{i+1}, \; c_i \to a_{i-1}$ in $\hat{I}_{p+r}$.

\medskip

Consider, for example, the case $n=7$. Then all admissible codes of weight $1$ are $142,124,1222,34,223$.
The corresponding sum is
$$
\hat{I}_1 = \sum_{\text{cyclic}} (-a_1 b_0-a_5 b_0+a_5 b_0 b_2 b_4-c_0+c_0 b_2 b_4).
$$
}
\end{remark}

\begin{proposition}\label{abc-int}
For odd $n$ one has
$$
I_0 = \prod_{j=0}^{n-1} a_j\,, \quad I_q = \hat{I}_1, \quad I_i = \hat{I}_{n-2i}\,,
$$
$$
G_0 = \prod_{j=0}^{n-1} c_j\,, \quad G_q = \hat{G}_1, \quad G_i = \hat{G}_{n-2i}\,.
$$
\end{proposition}
\proof
The proof is analogous to the proof of Proposition 5.3 in~\cite{OST99}.
\proofend

\medskip


\subsection{The spectral curve and invariant tori}

The spectral curve is a crucial component of algebraic-geometric integrability.
Below we always assume it to be generic. (As everywhere in this paper, ``generic'' means  values of parameters from some Zariski open subset in the space of parameters.)
It has a natural torus,   its Jacobian,  associated with it. It turns out that one can recover
a Lax function from the spectral curve and a point on the Jacobian, and vice versa: in our situation this correspondence is locally one-to-one.
The dynamics of the pentagram map becomes very simple on the Jacobian.
In this section, we construct this correspondence and describe the dynamics of the pentagram map.

\begin{definition}
{\rm
A {\it Floquet-Bloch solution} $\psi_{i,t}$ of a difference equation $\psi_{i+1,t} = \tilde{L}_{i,t} \psi_{i,t}$ is an eigenvector
of the monodromy operator:
$$
\tilde{T}_{i,t} \psi_{i,t} = k \psi_{i,t}.
$$
The normalization $\sum_{j=1}^4 \psi_{0,0,j} \equiv 1$ (i.e., the sum of all components of the vector $\psi_{0,0}$ is equal to $1$)
determines all vectors $\psi_{i,t}$ with $i,t \ge 0$ uniquely.
Denote the {\it normalized vectors} $\psi_{i,t}$ by $\Bar{\psi}_{i,t}$, i.e., $\Bar{\psi}_{i,t} = \psi_{i,t}/\left(\sum_{j=1}^4 \psi_{i,t,j}\right)$. (The vectors $\psi_{0,0}$ and $\Bar{\psi}_{0,0}$ are identical in this notation.) We also denote by $D_{i,t}$ the pole divisor of $\Bar{\psi}_{i,t}$ on $\Gamma$.
}
\end{definition}

\begin{remark}\label{Labc-xyz}
{\rm
We use the Lax function $\tilde{L}_{i,t}$ and the monodromy operator $\tilde{T}_{i,t}$ in the above definition to allow for both even and odd $n$. In the case of odd $n$ one can instead employ the Lax function  $L_{i,t}$ and the monodromy operator  $T_{i,t}$, while all the statements and proofs below  remain valid.
}
\end{remark}

\begin{theorem}\label{thm:floquet}
A Floquet-Bloch solution $\Bar{\psi}_{i,t}$ is a meromorphic vector function on $\Gamma$.
Generically its pole divisor $D_{i,t}$ has degree $g+3$.
\end{theorem}

\proof
The proof of the fact that the function $\Bar{\psi}_{i,t}$ is meromorphic on the spectral curve
$\Gamma$, as well as that its number of poles is $\deg D_{i,t}=\nu/2$, is identical to the proof of Proposition 3.4 in~\cite{FS11}.
The number $\nu$ of the branch points of $\Gamma$ is different:  in Theorem \ref{thm:genus}
we found that $2-2g=8-\nu$, where $g$ is the genus of the spectral curve.
This implies the required expression:  $\deg D_{i,t}=g+3$.
\proofend

\begin{definition}\label{def:spectral_map}
{\rm
Let $J(\Gamma)$ be the Jacobian of the spectral curve $\Gamma$, and $[D_{0,0}]$
is the equivalence class of the divisor $D_{0,0}$, the pole divisor of $\psi_{0,0}$, under the Abel map.
The pair consisting of the spectral curve $\Gamma$ (with marked points $O_i$ and $W_i$) and  a point $[D_{0,0}] \in J(\Gamma)$ is called the {\it spectral data.} The {\it spectral map} $S$   associates to a given generic twisted $n$-gon in $\CP^3$ its spectral data $(\Gamma, [D_{0,0}])$. 
}
\end{definition}

The   algebraic-geometric integrability is based on the following theorem.

\begin{theorem}\label{spectral-th1}
 For any $n$, the spectral map $S$ defines a bijection between 
a Zariski open subset of the space $\mathcal{P}_n = \{(x_i,y_i,z_i),\, 0 \le i \le n-1\}$
 and a Zariski open subset of the spectral data.
\end{theorem}
\begin{corollary}\label{cor-spec-th1}
 For odd $n$, the spectral map $S$ defines a bijection between a Zariski 
open subset of the space $\mathcal{P}_n \simeq \C^{3n}= \{(a_i,b_i,c_i),\, 0 \le i \le n-1\}$
 and a Zariski open subset of the spectral data.
\end{corollary}
\paragraph{Proof of Corollary~\ref{cor-spec-th1}.}
The statement follows from Theorem~\ref{spectral-th1}  and Definition~\ref{def:x-y-z} relating the coordinates $(x_i,y_i,z_i) $ and $(a_i,b_i,c_i)$ for odd $n$.
\proofend

The proof of Theorem~\ref{spectral-th1} is based on Proposition~\ref{prop-c} (which completes the construction of the direct spectral map) and Proposition~\ref{prop-b} (an independent construction of the inverse spectral map),
which are also used below to describe the corresponding pentagram dynamics.
It will be convenient to introduce the following notation for  divisors: $O_{pq}:=O_p+O_q$ and 
$W_{pq}:=W_p+W_q$ (e.g., $O_{12}:=O_1+O_2$).

 \begin{proposition}\label{prop-c}
The divisors of the coordinate functions $\psi_{i,t,1},..., \psi_{i,t,4}$ for $ 0 \le i \le n-1$ and any integer $t$ satisfy the following inequalities, provided that their divisors remain non-special up to time $t$:

For odd $n$ one has
  \begin{itemize}
  \item $(\psi_{i,t,1}) \ge -D +O_2 -i O_{23}+(i+1)W_{12}-t(W_{12}-O_{13});$
  \item $(\psi_{i,t,2}) \ge -D +(1-i)O_{23}+iW_{12}-t(W_{12}-O_{13});$
  \item $(\psi_{i,t,3}) \ge -D -i O_{23}+(i+1)W_{12}-t(W_{12}-O_{13});$
  \item $(\psi_{i,t,4}) \ge -D +O_2+(1-i) O_{23}+iW_{12}-t(W_{12}-O_{13});$
  \end{itemize}

For even $n$ one has
  \begin{itemize}
  \item $(\psi_{i,t,1}) \ge -D +O_2 + \lfloor\dfrac{i-t+2}{2} \rfloor W_{12}+ \lfloor\dfrac{i-t+1}{2} \rfloor W_{34}- \lfloor\dfrac{i+1}{2} \rfloor O_{24}-\lfloor\dfrac{i}{2} \rfloor O_{34}+tO_{14};$
  \item $(\psi_{i,t,2}) \ge -D + \lfloor\dfrac{i-t+1}{2} \rfloor W_{12}+ \lfloor\dfrac{i-t}{2} \rfloor W_{34}- \lfloor\dfrac{i-1}{2} \rfloor O_{24}-\lfloor\dfrac{i}{2} \rfloor O_{34}+tO_{14};$
  \item $(\psi_{i,t,3}) \ge -D + \lfloor\dfrac{i-t+2}{2} \rfloor W_{12}+ \lfloor\dfrac{i-t+1}{2} \rfloor W_{34}- \lfloor\dfrac{i+1}{2} \rfloor O_{34}- \lfloor\dfrac{i}{2} \rfloor O_{24}+tO_{14};$
  \item $(\psi_{i,t,4}) \ge -D +O_2 + \lfloor\dfrac{i-t+1}{2} \rfloor W_{12}+ \lfloor \dfrac{i-t}{2} \rfloor W_{34}- \lfloor\dfrac{i-1}{2} \rfloor O_{34}- \lfloor\dfrac{i}{2} \rfloor O_{24}+tO_{14};$
  \end{itemize}
  where $D=D_{0,0}$  corresponds to the divisor at $t=0$ and 
is an effective divisor of degree $g+3$, while $\lfloor x \rfloor$ is the floor (i.e., the greatest integer) function of $x$.

 \end{proposition}
 \proof
 The proof is a routine comparison of power expansions in $\lambda$ 
at the points $O_p, W_q$ for  
$k_i$ and $L_{i,t}$ and is very similar to the proof of Proposition 3.10 
in the 2D case in \cite{FS11}, although the 3D explicit expressions  are more involved. See more details in Appendix \ref{spec-data}.
 \proofend

\begin{proposition}\label{prop-b}
For any $n$,  given a generic spectral curve with marked points and a generic divisor $D$ of degree $g+3$
one can recover a sequence of matrices
\[
\tilde{L}_{i,t}(\lambda) =
\begin{pmatrix}
  0               & 0   & 0           & -1\\
  \lambda x_i y_i & 0   & 0           & 1\\
  0               & z_i & 0           & 1\\
  0               & 0   & \lambda x_i & 1
 \end{pmatrix}^{-1},
\]
for $ 0 \le i \le n-1 $ and any $ t.$
\end{proposition}

We describe  the reconstruction procedure and  prove this proposition in Appendix \ref{spec-data}.

\paragraph{Proof of Theorem~\ref{spectral-th1}.}
The proof consists of constructions of the spectral map $S$
and its inverse. The spectral map was  described in Definition \ref{def:spectral_map} based on Theorem \ref{thm:floquet}. 
We comment on an {\it independent} construction of the inverse spectral map now.

Pick an arbitrary divisor $D$ of degree $g+3$ in the equivalence class $[D_{0,0}] \in J(\Gamma)$ and apply
Proposition~\ref{prop-b}. ``A Zariski open subset of the spectral data'' is defined by spectral functions which may be singular
only at the points $O_i,W_i$ and by such divisors $[D] \equiv [D_{0,0}] \in J(\Gamma)$ that all divisors in
Proposition 6.17 with $0 \le i \le n-1$ up to time $t$ are non-special.
\proofend

The next theorem describes the time evolution of the pentagram map in the Jacobian of $\Gamma$.
The difference between even and odd $n$ is very similar to the 2-dimensional case.
Combined with Theorem~\ref{spectral-th1}, it proves the algebraic-geometric integrability of the 3D pentagram map.
(It also implies that it is possible to obtain explicit formulas of the coordinates of the pentagram map as functions of time using the Riemann $\theta$-functions.)

\begin{theorem}\label{time-evol}
 The equivalence class $[D_{i,t}]\in J(\Gamma)$ of the pole divisor $D_{i,t}$ of $\Bar{\psi}_{i,t}$ has the following time evolution:
 \begin{itemize}
  \item when $n$ is odd,
  $$
  [D_{i,t}] = [D_{0,0} - tO_{13} + iO_{23} + (t-i)W_{12}] ,
  $$
  \item when $n$ is even,
 $$
    [D_{i,t}] = \left[D_{0,0} - tO_{14} + \lfloor \dfrac{i+1}{2} \rfloor O_3 + \lfloor \dfrac{i}{2}\rfloor O_2 + i O_4 -\lfloor \dfrac{i-t+1}{2}\rfloor W_{12} - \lfloor  \dfrac{i-t}{2}\rfloor W_{34}\right].
$$
 \end{itemize}
 where $\deg{D_{i,t}} = g+3$ and $\lfloor x\rfloor $ is the floor function of $x$, and provided that
 spectral data remains generic up to time $t$.

 For an odd $n$ this discrete time evolution in $J(\Gamma)$ takes place along a straight line, whereas for an even $n$ the evolution
 goes along a ``staircase'' (i.e., its square goes along a straight line).
\end{theorem}
\proof
The vector functions $\psi_{i,t}$ with $i,t \ne 0$ are not normalized. The normalized vectors are equal to
$\Bar{\psi}_{i,t} = \psi_{i,t}/f_{i,t},$ where $f_{i,t}=\sum_{j=1}^4 \psi_{i,t,j}$.
Proposition~\ref{prop-c} implies that the divisor of each function $f_{i,t}$ is:
\begin{itemize}
\item for odd $n$,
$$
 (f_{i,t}) = D_{i,t} - D_{0,0} + tO_{13} - iO_{23}+ (i-t)W_{12},
$$
\item for even $n$,
$$ (f_{i,t}) = D_{i,t} - D_{0,0} + tO_{14} - \lfloor \dfrac{i+1}{2}\rfloor O_3 - \lfloor \dfrac{i}{2}\rfloor O_2 - i O_4+\lfloor  \dfrac{i-t+1}{2}\rfloor W_{12}+\lfloor \dfrac{i-t}{2}\rfloor W_{34}.
$$
\end{itemize}
Since the divisor of any meromorphic function is equivalent to zero, the result of the theorem follows.
The staircase dynamics is related to alternating jumps in the terms $\lfloor (i-t+1)/2\rfloor $ and $\lfloor (i-t)/2\rfloor $ as $t$ increases over integers.
\proofend

Note that although the pentagram map preserves the spectral curve, it exchanges the marked points.
The ``staircase'' dynamics on the Jacobian appears after the identification of curves with different marking. One cannot observe this dynamics in the  space of twisted polygons $\mathcal{P}_n$, before the application of the spectral map.

\section{Ramifications: closed polygons and symplectic leaves}

\subsection{Closed polygons}

Closed polygons in $\CP^3$ correspond to the monodromies $M=\pm\text{Id}$  in $SL(4,\C)$.
They form a subspace of codimension $15=\dim SL(4,\C)$ in the space of all twisted polygons $\mathcal{P}_n$.   The pentagram map on closed polygons in 3D is defined for $n \ge 7$.

\begin{theorem}\label{thm:closed}
  Closed polygons  in $\CP^3$ are singled out by the condition that either $(\lambda,k)=(1,1)$
  or $(\lambda,k)=(1,-1)$ is a quadruple point of $\Gamma$. Both conditions are equivalent to 9 independent linear constraints on $I_j,J_j,G_j$.
  Generically, the  genus of $\Gamma$ drops to $g=3q-9$ when $n$ is even, and to $g=3q-6$ when $n$ is odd, where $q=\lfloor n/2 \rfloor$.
  The dimension of the Jacobian $J(\Gamma)$ drops by $6$ for closed polygons for any $n$.
  Theorem~\ref{thm:floquet} holds with this genus adjustment, and
  Theorems~\ref{spectral-th1} and~\ref{time-evol} hold
  verbatim for closed polygons (i.e., on the subspace of closed polygons $\mathcal{C}_n \subset \mathcal{P}_n$).
 \end{theorem}

 \proof
For a twisted $n$-gon its monodromy matrix at a moment $t$ is equal to  $T_{0,t}(1)$ in the $(a,b,c)$-coordinates or to $\tilde{T}_{0,t}(1)$ in the $(x,y,z)$-coordinates.
An $n$-gon is closed if and only if $T_{0,t}(1)=\text{Id}$ or $T_{0,t}(1)=-\text{Id}$
(respectively, $\tilde{T}_{0,t}(1)=I_0 \,\text{Id}$ or $\tilde{T}_{0,t}(1)=-I_0 \,\text{Id}$).
For our definition of the spectral function, either of these conditions,  $T_{0,t}(1)=\pm\text{Id}$ or $\tilde{T}_{0,t}(1)=\pm I_0 \,\text{Id}$, implies that
$(\lambda,k)=(1,\pm 1)$ is a self-intersection point for $\Gamma$.

 The algebraic conditions implying that $(1,\pm 1)$ is a quadruple point are:
 \begin{itemize}
 \item $R(1,\pm 1) = 0$,
 \item $\partial_k R(1,\pm 1) = \partial_\lambda R(1,\pm 1) = 0$,
 \item $\partial_k^2 R(1,\pm 1) = \partial_\lambda^2 R(1,\pm 1) = \partial_{k \lambda}^2 R(1,\pm 1) = 0$,
 \item $\partial_k^3 R(1,\pm 1) = \partial_\lambda^3 R(1,\pm 1) = \partial_{kk \lambda}^3 R(1,\pm 1)= \partial_{k \lambda \lambda}^3 R(1,\pm 1) = 0$.
 \end{itemize}
 However, the function $R(\lambda,k)$ is special at the points $(1,\pm 1)$,
 because the following relation holds:
 $$
 R(1,\pm 1)=\pm\partial_k R(1,\pm 1) - \dfrac{1}{2} \partial_k^2 R(1,\pm 1) \pm \dfrac{1}{6} \partial_k^3 R(1,\pm 1).
 $$
 Consequently, the above 10 conditions are equivalent to only 9 independent linear equations on $I_j,J_j,G_j, 0 \le j \le q$.

 The proofs of Theorems~\ref{spectral-th1} and~\ref{time-evol} apply, mutatis mutandis, to
 the periodic case. To define the Zariski open set of spectral data for closed polygons, we confine to  spectral functions that can be singular only at the point $(\lambda,k)=(1,1)$ or $(1,-1)$
in addition to singularities at  $O_i$ and $W_i$ and use the same restrictions on divisors $D$ as in the proof of Theorem~\ref{spectral-th1}.

In  the periodic case we also have to  adjust
the count of the number $\nu$ of branch points  of $\Gamma$ and
 the corresponding calculation for the genus $g$ of $\Gamma$, cf. Theorem~\ref{thm:floquet}.
 Namely, as before, the function $\partial_k R(\lambda,k)$ has poles of total order $9n$ over $\lambda=0$, and zeroes of total order $6n$
 over $\lambda=\infty$. Now since $R(\lambda,k)$ has a quadruple point $(1,\pm 1)$, $\partial_k R(\lambda,k)$ has a triple zero at $(1,\pm 1)$.
 But $\lambda=1$ is not a branch point of $\Gamma$. Consequently, $\partial_k R(\lambda,k)$ has triple
 zeroes on 4 sheets of $\Gamma$ over $\lambda=1$.
 The Riemann-Hurwitz formula is $2-2g=8-\nu$, where the number of branch points
 for even $n$ is $\nu = 9n-6n-12=3n-12$, while for odd $n$ it is
 $\nu = 9n-6n-12+3=3n-9$. Therefore, we have $g=3q-9$ for even $n$, and $g=3q-6$ for odd $n$.
 \proofend


 \subsection{Invariant symplectic structure and symplectic leaves}\label{S:lax-sform}

 It was proved in~\cite{FS11} that in the 2D case an invariant symplectic structure
 on the space of twisted polygons ${\mathcal P}_n$ provided by Krichever-Phong's universal formula~\cite{KP97,KP98}
 coincides with the inverse of the invariant Poisson structure found in~\cite{OST99} when restricted to the symplectic leaves.
 We show that in 3D the same formula also provides an invariant symplectic structure defined on  leaves described below.
While we do not compute the symplectic structure explicitly in the coordinates
 $(a_i,b_i,c_i)$ or $(x_i,y_i,z_i)$  due to complexity
 of the formulas, the proofs are universal and  applicable in the higher-dimensional case of $\CP^d$ as well.  Finding an explicit expression of the symplectic structure or of the corresponding Poisson structure is still an open problem.

  \begin{definition}[\cite{KP97,KP98}]\label{KP-universal}
{\rm
 Krichever-Phong's universal formula defines a {\it pre-symplectic form} on the space of Lax operators, i.e.,
 on the space $\mathcal{P}_n$.
 It is given by the expression:
 $$
 \omega = -\dfrac{1}{2} \sum_{\lambda=0,\infty} {\text{res}} \thinspace \text{Tr}\left( \Psi_0^{-1} \tilde{T}_0^{-1} \delta \tilde{T}_0 \wedge \delta \Psi_0
 \right) \dfrac{d\lambda}{\lambda}.
 $$
 The matrix $\Psi_{0}:=\Psi_{0,t}(\lambda)$ is composed of the eigenvectors $\psi_{0,t}$ on different sheets of $\Gamma$ over the $\lambda$-plane,
 and it diagonalizes the monodromy matrix $\tilde{T}_{0}:=\tilde{T}_{0,t}(\lambda)$. (In this definition we drop the index $t$, because
 all variables correspond to the same moment of time.)

 The {\it leaves} of the 2-form $\omega$ are defined as submanifolds of $\mathcal{P}_n$, where the expression $\delta \ln{k}\, (d\lambda/\lambda)$
 is holomorphic. The latter expression is considered as a 1-form on the spectral curve $\Gamma$.
}
 \end{definition}

 \begin{proposition}\label{symp-leaves}
For even $n$  the leaves are singled out by $6$ conditions:
 $$
 \delta I_0 = \delta I_q = \delta G_0 = \delta G_q = \delta J_0 = \delta J_q = 0;
 $$
 For odd $n$ the leaves are singled out by $3$ conditions:
  $$
 \delta G_0 = \delta I_0 = \delta J_q = 0.
 $$
 \end{proposition}
 \proof
 These conditions follow immediately from the definition of the leaves and Lemma~\ref{spec-sing}.
 For example, at the point $O_1$ we have
 $$
 \delta \ln{k_1} \dfrac{d\lambda}{\lambda} = \left(\dfrac{1}{\lambda} \dfrac{\delta I_0}{I_0} + {\mathcal O}(1) \right) d\lambda.
 $$
 This 1-form is holomorphic in $\lambda$ if and only if $\delta I_0 = 0$.
 Similarly, we obtain $\delta(I_0/G_0)=0$ at the point $O_2$ for odd $n$. (One has to keep in mind that the local
 parameter around this point is $\lambda^{1/2}$.)
 \proofend

 \begin{remark}
 {\rm
 The definition of a presymplectic structure $\omega$ on $\mathcal{P}_n$ uses $\Psi_0$ and
 $\tilde{T}_0$ and hence relies on the normalization of $\Psi_0$.
 When restricted to the leaves from Proposition~\ref{symp-leaves}, the 2-form $\omega$ becomes
 independent of the normalization of the Floquet-Bloch solutions.
 Additionally, the form $\omega$ becomes non-degenerate, i.e., symplectic, when restricted to these leaves, as we prove below.
 The symplectic form is invariant with respect to the evolution given by the Lax equation.
 The proof is very similar to that of Corollary 4.2 in~\cite{IK02}
 (cf.~\cite{KP97,KP98} for other proofs).
 }
 \end{remark}

 \begin{theorem}\label{thm:rank}
  The rank of the invariant 2-form $\omega$  restricted to the leaves of Proposition~\ref{symp-leaves}
  is equal to $2g$.
 \end{theorem}
 \proof
 Since the 1-form $\delta \ln{k} \,d\lambda/\lambda$ is holomorphic on $\Gamma$, it can be represented as a sum of the basis
 holomorphic differentials:
\begin{equation}\label{kpformtmp4}
 \delta \ln{k} \,\frac{d\lambda}{\lambda} = \sum_{i=1}^{g} \delta U_i \,d\omega_i,
 \end{equation}
 where $g$ is the genus of $\Gamma$.
The coefficients $U_i$ can be found by integrating the last expression over the basis cycles $a_i$ of $H_1(\Gamma)$:
$$
U_i = \oint_{a_i} \ln{k} \,\frac{d\lambda}{\lambda}\, .
$$
According to formula~(5.7) in~\cite{KP00}, we have:
$$
\omega=\sum_{i=1}^{g+3} \delta \ln{k(p_i)} \wedge \delta \ln{\lambda(p_i)},
$$
where the points $p_i \in \Gamma,\; 1 \le i \le g+3,$ constitute the pole divisor $D_{0,0}$ of the
normalized Floquet-Bloch solution $\psi_{0,0}$.

After rearranging the terms, we obtain:
$$
\omega = \delta\left( \sum_{s=1}^{g+3} \int^{p_s} \delta \ln{k} \,\frac{d\lambda}{\lambda} \right)=
\delta\left( \sum_{s,i} \int^{p_s} \delta U_i \,d\omega_i \right) = \sum_{i=1}^{g} \delta U_i \wedge \delta\varphi_i,
$$
where
$$
\varphi_i = \sum_{s=1}^{g+3} \int^{p_s} d\omega_i
$$
are coordinates on the Jacobian  $J({\Gamma})$.
The variables $U_i$ and $\varphi_i$ are  natural Darboux coordinates for $\omega$, which also
turn out to be action-angle coordinates for the pentagram map. (The latter follows from the general properties of the Krichever-Phong universal form for a given Lax representation, cf. \cite{KP97,KP98}.)

Let us show that the functions $U_i$ are independent. Assume the contrary,
then there exists a vector $v$ on the space ${\mathcal P}_n$, such that $\delta U_i(v)=0$ for all $i$.
Then it follows from~(\ref{kpformtmp4}) that $\partial_v k \equiv 0$.
After applying the operator  $\partial_v$ to $R(\lambda,k)$, we conclude that $k$ satisfies
an algebraic equation of degree $3$, which is impossible, since $\Gamma$ is
a $4$-fold cover of the $\lambda$-plane.
 \proofend

\begin{remark}
{\rm
In more details, there are the following two cases:
 \begin{itemize}
 \item even $n=2q$. The dimension of the space $\mathcal{P}_n$ is $6q$. The codimension of the leaves is $6$.
 Therefore, the dimension of the leaves matches the doubled  dimension of the tori: $2g=6q-6$.
 \item odd $n=2q+1$. The dimension of the space $\mathcal{P}_n$ is $6q+3$. The codimension of the leaves is $3$.
 Again, the dimension of the leaves matches  the doubled  dimension of the tori: $2g=6q$.
 \end{itemize}
 }
 \end{remark}

The algebraic-geometric integrability in the complex case implies Arnold-Liouville integrability 
in the real one. Indeed, the pre-symplectic form depends on entries of the monodromy matrix 
in a rational way, since  it is independent of the permutation of sheets of 
the spectral curve $\Gamma$. Therefore, its restriction
to the space of the real $n$-gons provides a real pre-symplectic structure.
One obtains invariant Poisson brackets on the space of polygons $\mathcal{P}_n$
by inverting the real symplectic structure on the leaves, while employing invariants 
of Proposition \ref{symp-leaves} as the corresponding Casimirs.

\begin{problem}
Find an explicit formula for an invariant Poisson structure with the above symplectic leaves.
\end{problem}

 \medskip


\section{A Lax representation in higher dimensions}\label{S:higher-lax}

The origin of the integrability of the pentagram map is the presence
of its scaling invariance.  Assume that $gcd(n,d+1)=1$.
The difference equation (\ref{eq:difference_anyD})
$$
V_{j+d+1} = a_{j,d} V_{j+d} + a_{j,d-1} V_{j+d-1} +...+ a_{j,1} V_{j+1} + (-1)^d V_j
$$
allows one to introduce coordinates $a_{j,1},a_{j,2},...a_{j,d}, \, 0\le j\le n-1$, on the space
of twisted $n$-gons in any dimension $d$.

\begin{proposition-conjecture} {\bf (The scaling invariance)}
The pentagram map on twisted $n$-gons in $\CP^d$ is invariant with respect to the following scaling transformations:
\begin{itemize}
\item for odd $d=2\varkappa+1$ the transformations are
$$
a_{j,1} \to s a_{j,1},\; a_{j,3} \to s a_{j,3},\;a_{j,5} \to s a_{j,5},\;...\;, a_{j,d} \to s a_{j,d}\,,
$$
while other coefficients $a_{j,2l}$ with $l=1,...,\varkappa$ do not change;

\item for even $d=2\varkappa$ the transformations are
$$
a_{j,1} \to s^{-\varkappa}a_{j,1} ,\; a_{j,3} \to
s^{1-\varkappa}a_{j,3}, \;...\:, a_{j,d-1} \to  s^{-1}a_{j,d-1},
$$
$$
a_{j,2} \to sa_{j,2} ,\; a_{j,4} \to s^{2}a_{j,4},
\;...\:,  a_{j,d} \to s^{\varkappa}a_{j,d}
$$
\end{itemize}
for all $s\in \C^*$.\footnote{We thank G.Mari-Beffa for correcting an
error in the scaling for even $d$ in the first version of this manuscript, as well
as in the short version \cite{16}. This error related to numerics with
a different choice of vertices for the diagonal planes leads to another 
system, different from $T_{p,r}$, which also turns out to be integrable and will be discussed elsewhere.}
\end{proposition-conjecture}

\proof
In any dimension $d$ the pentagram map is a composition of involutions $\alpha$ and $\beta$, see Section \ref{sect:invol}.
(More precisely, $\alpha$ is not an involution for even $d$, but its square $\alpha^2$ is a  shift in the vertex index, see \cite{OST99} for the 2D case.)
One can prove  that the involution $\alpha: V_i \to W_i $  in any dimension has the form
$$
W_{j+d+1} = (-1)^{d+1}(a_{\star,1} W_{j+d} + a_{\star,2} W_{j+d-1} +...+ a_{\star,d} W_{j+1} - W_j)\,,
$$
where $\star$ stands for the first index, which is irrelevant for the scaling (Lemma \ref{map-alpha} proves the case $d=3$).

We call this Proposition-conjecture because the proof of an analog of  Lemma \ref{map-beta} (for the map $\beta$) in higher dimensions is computer assisted.
One verifies that for a given dimension $d$ the coefficients consist of the terms that are consistent with the scaling. 
\proofend

We obtained explicit formulas, and hence a  direct (theoretical) proof of the scaling invariance for the pentagram maps up to dimension $d\le 6$. This bound is related to computing powers to produce explicit formulas and might be extended. However, we have no general purely theoretical proof valid for all $d$ and it would be very interesting to find it.

\begin{problem}
Find a general proof of the scaling invariance of the pentagram map in any dimension $d$.
\end{problem}

\begin{theorem}\label{thm:lax_anyD}
The scale-invariant pentagram map on twisted $n$-gons in any dimension $d$ is a completely integrable system. It is described by the  Lax matrix
\[
L_j^{-1}(\lambda) =
\left(
\begin{array}{cccc|c}
0 & 0 & \cdots & 0    &(-1)^d\\ \cline{1-5}
\multicolumn{4}{c|}{\multirow{4}*{$D(\lambda)$}} & a_{j,1}\\
&&&& a_{j,2}\\
&&&& \cdots\\
&&&& a_{j,d}\\
\end{array}
\right),
\] 
where $D(\lambda)$ is the following diagonal matrix of size $d \times d$:
\begin{itemize}
\item for odd $d=2\varkappa+1$, one has
$D(\lambda) = \text{diag}(\lambda,1,\lambda,1,...,\lambda)$;

\item for even $d=2\varkappa$, one has
one has $D(\lambda) = \text{diag}(1,\lambda,1,\lambda,...,1, \lambda)$.
\end{itemize}
\end{theorem}

\proofsketch
By using the scaling invariance of the pentagram map, one derives the Lax matrix  exactly
in the same way as in 3D, see Section \ref{S:lax-func}.
Namely, first  construct the $(d+1) \times (d+1)$-matrix $N_j(s)$
depending on our scaling parameter $s$,
and then use the formula $L_j^{-1}(\lambda) = \left(g^{-1} N_j(s) g \right)/s^m$ with a suitable choice
of the diagonal $(d+1) \times (d+1)$-matrix $g$ and an appropriate function of the parameter $s$.

For odd $d=2\varkappa+1$, we have $g = \text{diag}(1,s,1,s,...,1,s)$, $m=1$, and $\lambda \equiv s^{-2}$,  whereas
for even $d=2\varkappa$, we have $g = \text{diag}(1, s^{-\varkappa},s,s^{1-\varkappa},s^2,
..., s^{\varkappa-1},s^{-1},s^{\varkappa})$,  $m={\varkappa}$, and $\lambda \equiv s^{-d-1}$.
The  Lax representation with a spectral parameter is constructed as we described above.

Using the genericity assumptions similar to those used in the 2D and 3D cases, one constructs the spectral map and its inverse, which is equivalent to algebraic-geometric integrability of the pentagram map. Coefficients of the spectral curve form a maximal family of first integrals. Along with a (pre)symplectic structure defined by the Krichever--Phong formula, this provides the Arnold--Liouville integrability  of the system on the corresponding symplectic leaves in the real case.
\proofend

\medskip
The scaling parameter has a clear meaning in the continuous limit:
\begin{proposition}
For any dimension $d$ the continuous limit of the scaling transformations
corresponds to the spectral shift $L\rightarrow L+\lambda$ of the differential operator $L$.
\end{proposition}

\proof
In 2D this was proved in \cite{OST99}. A continuous analog of the difference equation (\ref{eq:difference_anyD}) is
$$
G(x+(d+1)\ep)=a_d(x,\ep)G(x+d\ep)+...+a_1(x,\ep)G(x+\ep)+(-1)^d G(x),
$$
where $G(x)$ satisfies the differential equation (\ref{eq:diff_anyD_onG}) with a differential operator $L$ of the form (\ref{operatorL}).
Using the Taylor expansion for $G(x+j\ep)$ and the expansion $a_j(x,\ep)=a_j^0(x)+\ep a_j^1(x)+...,$
we obtain expressions of $a_k^i$ in terms of the coefficients of $L$, i.e., in terms of functions
$u_j(x)$ and their derivatives.
We find that the terms $a_k^0$ are constant, $a_k^1=0$ for all $k$, while $a_k^i$ for $i\ge 2$ are linear in $u_{d-i+1}$ and  differential polynomials
in the preceding coefficients $u_{d-1},...,u_{d-i+2}$.

The scaling parameter also has an expansion in $\ep$:  $s=\tau_0+\ep \tau_1+\ep^2 \tau_2+...$.
We apply it to the coefficients $a_i(x,\ep)$ and impose
the condition that $a_k^0$ and $a_k^1=0$ are fixed, similarly to \cite{OST99}.
By term-wise calculations (different in the cases of even and odd $d$ and using the ``triangular form" of the expressions for  $a_k^i$),
one successively obtains that $\tau_0=1$,\, $\tau_1=...=\tau_d=0$, i.e.,
$s$ must have the form $s=1+\tau_{d+1} \ep^{d+1}+{\mathcal O}(\ep^{d+2})$.
Its action shifts only the last term of $L$: $u_0\to u_0+{\rm const}\cdot \tau_{d+1} \ep^{d+1}$, i.e.,
it is equivalent to the spectral shift $L\rightarrow L+\lambda$.
\proofend

\smallskip

Note that the spectral shift commutes with the KdV flows. Indeed,
$d/dt (L+\lambda)=d/dt\, L=[Q_2,L]=[Q_2, L+\lambda]$, since
$Q_2(L):=\partial^2 +\frac{2}{d+1}u_{d-1}=Q_2(L+\lambda)$ for operators $L$ of degree $d+1\ge3$.
Equivalently, the pentagram map commutes with the scaling transformations in the continuous limit.

\medskip



 \section{Appendices}

\subsection{Continuous limit in the 3D case}

In this section we present explicit formulas manifesting Theorem \ref{thm:cont-lim} on the continuous limit of the 3D pentagram map.
Consider a curve $G(x)$ in $\R^4$ given by the differential equation
$$
G''''+u(x)G''+v(x)G'+w(x)G=0
$$
with periodic coefficients $u(x),v(x),w(x)$.
To find the continuous limit, we fix $\epsilon$ and consider a plane $P_\ep(x)$ passing through the three points $G(x-\epsilon), G(x), G(x+\epsilon)$ on this curve.
We are looking for an equation of  the envelope curve $L_\epsilon(x)$ for these planes.

This envelope curve  $L_\ep(x)$ satisfies the following system of equations:
\begin{align*}
&\det | G(x), G(x+\ep), G(x-\ep),  L_\ep(x) |=0\\
&\det | G(x), G(x+\ep), G(x-\ep),  L'_\ep(x) |=0\\
&\det | G(x), G(x+\ep), G(x-\ep),  L''_\ep(x) |=0\,.
\end{align*}
By considering the Taylor expansion and using the normalizations
$\det | L_\ep, L'_\ep, L''_\ep, L'''_\ep |=1$ and $\det | G, G', G'', G''' |=1$
we find that
\begin{equation}\label{3D-L}
L_\ep(x) = G(x)+ \dfrac{\ep^2}{6} \left( G''(x)+\dfrac{u}{2} G(x)
\right) +{\mathcal O}(\epsilon^4)
\end{equation}
as $\epsilon\to 0$.
Now, the equation $L''''_\ep+u_\ep L''_\ep+v_\ep L'_\ep+w_\ep L_\ep=0$ implies that:
\begin{align*}
u_\ep &= u + \dfrac{\ep^2}{3}(v'-u'')+{\mathcal O}(\ep^4),\\
v_\ep &= v + \dfrac{\ep^2}{6}(2w'+v''-uu'-2u''')+{\mathcal O}(\ep^4),\\
w_\ep &= w + \dfrac{\ep^2}{12}(2w''-uu''-vu'-u'''')+{\mathcal O}(\ep^4).
\end{align*}

These equations describe the $(2,4)$-equation in the Gelfand-Dickey hierarchy:

\[
\dot{L} = [Q_2,L] \Leftrightarrow
\begin{cases}
\dot{u} = 2v'-2u'',\\
\dot{v} = v''+2w'-uu'-2u''',\\
\dot{w} = w'' - \dfrac{1}{2} vu'-\dfrac{1}{2}uu'' - \dfrac{1}{2} u'''',
\end{cases}
\]
where
$L = \partial^4 + u \partial^2 + v \partial + w$ and $Q_2 = (L^{2/4})_+ = \partial^2 + \dfrac{u}{2}$.

\begin{remark}\label{other-choice}
{\rm
A different choice of the points defining the plane $P_\ep(x)$
on the original curve leads to the same continuous limit.
For instance, the choice of $G(x-3\ep),G(x+\ep),G(x+2\ep)$ results in the
same expression for $L_\ep(x)$, where in (\ref{3D-L}) instead of the coefficient
$\ep^2/6$ one has $7\ep^2/6$.
This leads to the  same evolution of the curve $G$ with a different time parameterization, cf. Remark \ref{rm:diffchoice}.
}
\end{remark}

\medskip


\subsection{Higher terms of the continuous limit}\label{cont2D}

 Recall that in the continuous limit for the pentagram map in $\RP^d$ the envelope for osculating planes moves according to the $(2,d+1)$-KdV equation
 (Theorem \ref{thm:kdv}).
 This evolution is  defined by the  $\ep^2$-term of the expansion of the function $L_\ep(x)$.

 The same proof works in the following more general setting. Let $L$ be a differential operator (\ref{operatorL})
 of order $d+1$ and $G$ a non-degenerate curve defined by its solutions: $LG=0$.

\begin{proposition}
Assume that the  curve $G$ evolves according to the law
 $ \dot G= Q_m G$, where $Q_m:=(L^{m/(d+1)})_+$ is the differential part of the $m${\rm th} power of the operator $Q=L^{1/(d+1)}$. Then this evolution defines the equation $\dot  L=[Q_m, L]$, which is the
 $(m,d+1)$-equation in the corresponding KdV hierarchy of $L$.
\end{proposition}

Furthermore, one can define the simultaneous evolution of all terms in the $\ep$-expansion of $L_\ep(x)$ using the following construction. For the pseudodifferential operator  $Q:=L^{1/(d+1)}$ consider the formal series
$\exp (\ep Q):=1+\ep Q +\frac{\ep^2}{2} Q^2+...$ and take its differential part:
$$
\left(\exp (\ep Q)\right)_+=\left(1+\ep Q +\frac{\ep^2}{2} Q^2+...\right)_+=1+\ep Q_1 +\frac{\ep^2}{2} Q_2+...=\sum_0^\infty \frac{\ep^m}{m!}Q_m\,.
$$
For each power of $\ep$ this is a multiple of the differential operator $Q_m$, which is the differential part of the $m${\rm th} power $Q^m$ of the operator $Q=L^{1/(d+1)}$.

\begin{corollary}
The formal evolution equation $ \dot G= \left(\exp (\ep Q)\right)_+ G$  corresponds to the full
KdV hierarchy $\dot  L=[\left(\exp (\ep Q)\right)_+ , L]$, where the operator $L$ is of order $d+1$
and the $(m,d+1)$-equation corresponds to the power  $\ep^m$.
\end{corollary}

A natural question is which equations of this hierarchy actually appear as the evolution of the envelope $L_\ep(x)$.
Recall that only even powers of $\ep$ arise in the  expansion  of the function $L_\ep(x)$ for the continuous limit of the
pentagram map. The $\ep^2$-term gives the $(2,d+1)$-KdV equation.
It turns out that the $\ep^4$-term in the continuous limit of the 2D pentagram map results in the equation
very similar to the $(4,3)$-equation in the KdV hierarchy (which is a higher-order Boussinesq equation).
Although the numerical coefficients in these differential equations are different, one may hope to
obtain the exact equations of the KdV hierarchy for different $m$ by using an appropriate rescaling.
This allows one to formulate

\begin{problem}
{\rm
Do higher $(m, d+1)$-KdV flows appear as the $\ep^{m}$-terms in the expansion of the envelope $L_\ep(x)$ for the continuous limit
of the pentagram map  for any even $m>2$?
}
\end{problem}


\subsection{Bijection of the  spectral map}\label{spec-data}

In this appendix we sketch the proof of Proposition \ref{prop-c} and
prove Proposition  \ref{prop-b}, which allows one to reconstruct the $L$-matrix from spectral data, and hence complete the proof of Theorem \ref{spectral-th1} on the spectral map.

\begin{proposition}{\bf (= Proposition \ref{prop-b})}
For any $n$,  given a generic spectral curve with marked points and a generic divisor $D$  of degree $g+3$
one can recover a sequence of matrices
\[
\tilde{L}_{i,t}(\lambda) =
\begin{pmatrix}
  0               & 0   & 0           & -1\\
  \lambda x_i y_i & 0   & 0           & 1\\
  0               & z_i & 0           & 1\\
  0               & 0   & \lambda x_i & 1
 \end{pmatrix}^{-1},
\]
for $ 0 \le i \le n-1 $ and any $ t.$
\end{proposition}

\proof
Without loss of generality we describe  the procedure to reconstruct the matrices $L_i(\lambda):=\tilde{L}_{i,0}(\lambda), $ for $ 0 \le i \le n-1$ and $t=0$.
  \begin{enumerate}
  \item First, we pick functions $\psi_{i,j}:=\psi_{i,0,j}$ for $\; 0 \le i \le n-1,\; 1 \le j \le 4,$ and $t=0$, satisfying Proposition~\ref{prop-c}.
  Note that according to the Riemann-Roch theorem, the functions $\psi_{i,1}$ and $\psi_{i,4}$ are defined up to
  a multiplication by constants, whereas the functions $\psi_{i,2}$ and $\psi_{i,3}$
  belong to 2-dimensional subspaces. The functions $\psi_{i,1}$ and $\psi_{i,4}$ belong to the same subspaces.
  We pick the pairs of functions $\psi_{i,1},\; \psi_{i,3}$ and $\psi_{i,2},\; \psi_{i,4}$ to be linearly independent. Observe that any sets of functions $\psi_{i,1},...,\psi_{i,4} $ satisfying Proposition~\ref{prop-c} are related
  by gauge transformations $\psi_i \to g_i^{-1} \psi_i$, where
  \[
  g_i =
  \begin{pmatrix}
   A_i & 0   & 0   & 0\\
   0   & B_i & 0   & E_i\\
   F_i & 0   & C_i & 0\\
   0   & 0   & 0   & D_i
  \end{pmatrix},\;
  g_i = g_{i+n},
  \]
  and $\psi_i$ stands for $\psi_i=(\psi_{i,1},...,\psi_{i,4})^T$. We also 
  define $\psi_n$ to be $\psi_n=I_0k\psi_0$ for any $n$ in 
  $(x,y,z)$-variables.
 
  \item We find the unique matrix $L'_i$ satisfying the equation $\psi_i = (L'_i)^{-1} \psi_{i+1}$:
  \[
  L'_i(\lambda) =
  \begin{pmatrix}
   0               & 0       & 0               & t_{i,1}\\
   \lambda t_{i,5} & 0       & \lambda t_{i,6} & t_{i,2}\\
   0               & t_{i,7} & 0               & t_{i,3}\\
   \lambda t_{i,8} & 0       & \lambda t_{i,9} & t_{i,4}
  \end{pmatrix}^{-1}.
\]

 \item One can check that there exists the unique choice of the matrices $g_i, \; 0 \le i \le n-1,$
  such that the equality $L_i(\lambda) = g_{i+1} L'_i(\lambda) g_i^{-1}$ is possible.
  The latter is equivalent to the following system of equations ($0 \le i \le n-1$):
  $$
  \dfrac{A_i t_{i,1}}{D_{i+1}}=-1; \quad \dfrac{B_i t_{i,2}+E_i t_{i,4}}{D_{i+1}}=1; \quad \dfrac{D_i t_{i,4}}{D_{i+1}}=1;
  $$
  $$
  \dfrac{B_{i+1} F_i t_{i,1}+B_{i+1} C_i t_{i,3}-C_i E_{i+1} t_{i,7}}{B_{i+1} D_{i+1}}=1; \quad B_i t_{i,6} + E_i t_{i,9}=0; \quad C_{i+1} t_{i,8} - F_{i+1} t_{i,9}=0.
  $$
  These equations decouple and may be solved explicitly. One only needs to check the solvability of $n$ equations
  $\dfrac{D_i t_{i,4}}{D_{i+1}}=1, \; 0 \le i \le n-1,$ for $n$ variables $D_i, \; 0 \le i \le n-1$.
  A non-trivial solution exists provided that $\prod_{i=0}^{n-1} t_{i,4} = 1$. It depends on an arbitrary constant, which corresponds to  multiplication of all matrices $g_i$ by the same number and does not affect the Lax matrices.
  One can check that
  $$
  t_{i,4} = \dfrac{\psi_{i+1,4}(O_1)} {\psi_{i,4}(O_1)} \text{ and } \prod_{i=0}^{n-1} t_{i,4} = \dfrac{\psi_{n,4}(O_1)} {\psi_{0,4}(O_1)}=I_0 k(O_1).
  $$
By using Lemma \ref{spec-sing} we find the value $k(O_1)  ={1}/{I_0}$ as required. Now the remaining variables $A_i,B_i,C_i,E_i,F_i,\; 0 \le i \le n-1,$ are uniquely determined.
 \end{enumerate}
\proofend

 \begin{corollary}\label{cor-d}
For odd $n$,  given a generic spectral curve with marked points and a generic divisor $D$
one can recover a sequence of matrices
\[
L_{i,t}(\lambda) =
 \begin{pmatrix}
  0       & 0   & 0       & -1\\
  \lambda & 0   & 0       & c_j\\
  0       & 1   & 0       & b_j\\
  0       & 0   & \lambda & a_j
 \end{pmatrix}^{-1}
\]
with  $0 \le i \le n-1$ and any $t$.
\end{corollary}

\proof
The statement follows from Proposition~\ref{prop-b} and the fact that
$(a_j,b_j,c_j),\; 0 \le i \le n-1,$ are coordinates on the space $\mathcal{P}_n$ for odd $n$.
\proofend

We complete the exposition with a sketch of the proof for Proposition \ref{prop-c} for even $n$
(the case of odd $n$ is similar).

\begin{proposition}{\bf (= Proposition \ref{prop-c}$'$)}
For even $n$, the divisors of the coordinate functions $\psi_{i,t,1},...,
\psi_{i,t,4}$ for $ 0 \le i \le n-1$ and any integer $t$ satisfy the
following inequalities,  provided that their divisors remain non-special up to time $t$:
  \begin{itemize}
  \item $(\psi_{i,t,1}) \ge -D +O_2 + \lfloor\dfrac{i-t+2}{2} \rfloor W_{12}+ \lfloor\dfrac{i-t+1}{2} \rfloor W_{34}- \lfloor\dfrac{i+1}{2} \rfloor O_{24}-\lfloor\dfrac{i}{2} \rfloor O_{34}+tO_{14};$
  \item $(\psi_{i,t,2}) \ge -D + \lfloor\dfrac{i-t+1}{2} \rfloor W_{12}+ \lfloor\dfrac{i-t}{2} \rfloor W_{34}- \lfloor\dfrac{i-1}{2} \rfloor O_{24}-\lfloor\dfrac{i}{2} \rfloor O_{34}+tO_{14};$
  \item $(\psi_{i,t,3}) \ge -D + \lfloor\dfrac{i-t+2}{2} \rfloor W_{12}+ \lfloor\dfrac{i-t+1}{2} \rfloor W_{34}- \lfloor\dfrac{i+1}{2} \rfloor O_{34}- \lfloor\dfrac{i}{2} \rfloor O_{24}+tO_{14};$
  \item $(\psi_{i,t,4}) \ge -D +O_2 + \lfloor\dfrac{i-t+1}{2} \rfloor W_{12}+ \lfloor \dfrac{i-t}{2} \rfloor W_{34}- \lfloor\dfrac{i-1}{2} \rfloor O_{34}- \lfloor\dfrac{i}{2} \rfloor O_{24}+tO_{14};$
  \end{itemize}
  where $D$ is an effective divisor of degree $g+3$, and $\lfloor x \rfloor$ is the floor  function of $x$.
\end{proposition}

\proof
First, we prove these inequalities for $t=0$ and $0 \le i \le n-1$.
For illustration we find the multiplicities of the components of the vector $\psi_{i,0}$ at the point $O_2$, while other points can be treated in a similar fashion.
We employ the matrices $\tilde{L}_{i,t}$ in the coordinates $x_i,y_i,z_i$.

Notice that a cyclic permutation of indices $(n-1,n-2,...,1,0)$ changes the monodromies 
$T_i \to T_{i+1}$ and the Floquet-Bloch solutions
$\bar{\psi}_i \to \bar{\psi}_{i+1}$. For even $n$, it also permutes
$\bar{\psi}_i(O_2) \leftrightarrow \bar{\psi}_i(O_3)$ and $W_{12} \leftrightarrow W_{34}$, i.e.,
the corresponding pairs of the vectors $\bar{\psi}_i$ at the points $(W_1,W_2)$ and $(W_3,W_4)$ are swapped.

Using the asymptotic expansion of $\tilde{T}_{0,t}(\lambda)$ at $\lambda=0$, 
the definition of the Floquet-Bloch solution,
and the normalization condition, one can show that $\psi_{0,0}=(O(\lambda),O(\lambda),1+O(\lambda),O(\lambda))^T$ as $\lambda\to 0$ at the point $O_2$.
Since 
$$
L_{1,0}(\lambda)L_{0,0}(\lambda) =
\begin{pmatrix}
  1   & 1   & 0 & 0\\
  0   & 0   & 0 & 0\\
  y_1 & y_1 & 0 & 0\\
  0   & 0   & 0 & 0
\end{pmatrix}
\dfrac{1}{x_0 x_1 y_0 y_1 \lambda^2} + O\left( \dfrac{1}{\lambda} \right) \text{ as } \lambda\to 0,
$$
and $\psi_{2,0}=L_{1,0} L_{0,0} \psi_{0,0}$, generically one has $\psi_{2,0}=(O(1),O(1),O(1/\lambda),O(1))^T$ at $O_2$.

By definition, the normalized vectors are $\bar{\psi}_{i,t} = f_{i,t} \psi_{i,t}$.
Using a cyclic permutation, we find that $\bar{\psi}_{2k,0}=(O(\lambda),O(\lambda),1+O(\lambda),O(\lambda))^T$
and that $f_{2,0}(\lambda) = O(\lambda)$ at $O_2$. Using the permutation argument again,
we derive that $f_{i+2,0}(\lambda)/f_{i,0}(\lambda) = O(\lambda)$ at $O_2$ for even $i$.
Therefore, one has $f_{2k,0}(\lambda) = O(\lambda^k)$ at $O_2$. Now the required multiplicities for the vector $\psi_{2k,0}$ at $O_2$ follow.
Furthermore, since $\psi_{2k+1,0} = L_{2k} \psi_{2k,0}$, one can check that generically $f_{2k+1,0}(\lambda)/f_{2k,0}(\lambda) = O(1)$
and $\bar{\psi}_{2k+1,0}=(O(1),O(\lambda),O(1),O(1))^T$ at the point $O_2$. This establishes also
the multiplicities for the vector $\psi_{2k+1,0}$ at $O_2$.

Having proved the proposition for $t=0$, one can prove it for $t>0$ by using the formula $\psi_{i,t+1} = \tilde{P}_{i,t} \psi_{i,t}$.
Note that it suffices to study the cases $t=0$ and $t=1$ only.
Consider, for example, the multiplicity of the function $\psi_{i,1,1}$ at the point $O_2$.
Since $\psi_{i,1,1}=(\psi_{i,0,2}-\psi_{i,0,4})/(x_i(1+y_i+z_{i+1}))$, one can check that the multiplicity of the right-hand side at $O_2$
is $1-k$ for $i=2k$ and it is equal to $-k$ for $i=2k+1$, i.e., $\psi_{i,1,1}$ and $\psi_{i,0,1}$ have the same multiplicities at $O_2$.
Other cases are treated in a similar way.
\proofend



\end{document}